%% file: SO_Bolten_Gottschalk_Hahn-Saadi.tex
\documentclass[12pt,twoside]{article}
\usepackage[utf8]{inputenc}
\usepackage[english]{babel}
\usepackage{amsmath}
\usepackage{amsfonts}
\usepackage{amssymb}
\usepackage{amsthm}
\usepackage{algorithm2e}
\usepackage{bbm}
\usepackage{graphicx}
\usepackage{lmodern}
\usepackage{todonotes}
\DeclareMathAlphabet{\mathpzc}{OT1}{pzc}{m}{it}
\usepackage{fourier}
\usepackage {subfig}
\usepackage[section]{placeins}
\usepackage[onehalfspacing]{setspace} 
\usepackage{booktabs}
\usepackage{wrapfig}
\usepackage{url}
\usepackage{pdfpages}
\usepackage{floatflt}
\usepackage{qtree}
\usepackage[shortlabels]{enumitem}
\usepackage{listings}

\usepackage[autooneside=false]{scrlayer-scrpage}
\usepackage{tikz}
\usetikzlibrary{trees}

\ihead{\headmark}
\ohead{\pagemark}
\theoremstyle{definition}

\theoremstyle{plain}

\newcommand{\spur}{\text{tr }}

\newcommand{\parableit}[2]{\frac{\partial {#1}}{\partial {#2}}}
\newcommand{\intOm}[1][ ]{\int_{\Omega}{#1}\,dx}

\numberwithin{equation}{section} 

\title{\textsc{Numerical shape optimization to decrease failure probability of ceramic
structures}}
\author{\textsc{M. Bolten, H. Gottschalk, C. Hahn and M. Saadi}\\
School of Mathematics and Natural Science\\ University of Wuppertal\thanks{\texttt{$\{$bolten@math.,hanno.gottschalk@,chahn@,saadi@$\}$uni-wuppertal.de}}}

\begin{document}

\maketitle
\begin{abstract}
 Ceramic is a material frequently used in industry because of its favorable properties. Common approaches in shape optimization for ceramic structures aim to minimize the tensile stress acting on the component, as it is the main driver for failure. In contrast to this, we follow a more natural approach by minimizing the component's probability of failure under a given tensile load. Since the fundamental work of Weibull, the probabilistic description of the strength of ceramics is standard and has been widely applied. Here, for the first time, the resulting failure probabilities are used as objective functions in PDE constrained shape optimization.  
 To minimize the probability of failure, we choose a gradient based method combined with a first discretize then optimize approach. For discretization finite elements are used. Using the Lagrangian formalism, the shape gradient via the adjoint equation is calculated at low computational cost. The implementation is verified by comparison of it with a finite difference method applied to a minimal 2d example. Furthermore, we construct shape flows towards an optimal / improved shape in the case of a simple beam and a bended joint. 
\end{abstract}



\section{Introduction}
\label{intro}
In his fundamental 1939 paper \cite{weibull}, W. Weibull summarized the experimental situation found in measurement of the ultimate strength of brittle materials: \emph{"The classical theory is obviously incompatible with numerous results of experimental research. This discrepancy may be bridged over by considering as an essential element of the problem the dispersion obtained in experimental measuring of the ultimate tensile strength (UTS). Viewed from this standpoint, the UTS of a material can not be expressed by a single numerical value, as has been tone heretofore, and a statistical distribution will be required for this purpose"}.   Since this time, the Weibull distribution has become an indispensable tool in the classification  of ceramics. In particular, the Weibull module $m$ which controls the dispersion of the UTS, has become a standard material property,  see e.g. \cite{munz}.   

In mechanical design, the ultimate tensile strength of a component has to be well above the loads applied. But even in situations, when the loads can be foreseen, Weibull's insight forces engineers to work with levels of reliability rather than with ultimate safety. When choosing shape and dimensions of an engineered part, the quest rather is to control a probability distribution than to keep the maximal stress below a given threshold. This logic has been successfully applied in technological areas ranging from space shuttle heat shields \cite{nasa}, gas turbine combustion chambers \cite{riesch} to dental prostheses \cite{dental}. This is particularly simple in the case of the Weibull model, since the failure distribution of the engineered component over the applied loads is modeled by scale and shape parameters, where the shape parameter is Weibull's module $m$ and only the scale quantity depends on the component's design. Maximizing the scale therefore corresponds to maximizing the probabilistic endurance or the component over the entire range of loads.  

The maximization of the Weibull scale of a component as a functional of the component's shape puts the control of failure for ceramic components in the framework of shape optimization \cite{allaire,bucur,chenais,haslinger,sokolowski}. In fact, as has been observed recently, the so obtained shape optimization problem has the favorable property that the objective functional is differentiable \cite{bolten,gottschalk,schmitz}. This is in sharp contrast to the peak stress criterion commonly used in optimization, which is non differentiable as objective functional due to  maximizing stress over all locations on the component. Gradient based shape optimization techniques, in particular in conjunction with the highly efficient calculation of shape derivatives via adjoint equations, have proven their potential in many engineering applications. In particular this applies to aero design, where the objective functionals have always been differentiable, see e.g. \cite{mohammadi,Frey2009}.   It is therefore natural to use the smoothing nature of probabilistic models to extend gradient based optimization to the design objective of (probabilistic) mechanical integrity as well.        

In the mathematical literature, shape and topology optimization for the linear elasticity PDE have been predominantly applied to the compliance functional, see e.g. \cite{allaire,bendsoe,conti}. This is however not directly related to mechanical integrity.

In the present article we for the first time minimize failure probabilities numerically using a first discretize, then adjoin strategy and apply it to a simple 2D design problem. We demonstrate that the  shape derivatives via the adjoint method can be calculated numerically with a  high level of precision. The resulting geometry flows, under volume constraints, are stable over rather significant changes of the geometry and converge to an (numerically) optimal solution. Remarkably this is true without any artificial smoothing of the shape gradients.  

The paper is organized as follows:  In Section \ref{FailProbs} we derive failure probabilities according from fundamental assumptions of elastic fracture mechanics and an initial flaw distribution using the a point process model introduced in \cite{bolten}. In this way one obtains Weibull's classical model with slight modifications \cite{batendorf, weibull}. For a derivation of the same model via extremal value statistics see e.g. \cite{riesch} and references therein.   We discuss the numerical approximation of the objective functional via the discretization of the PDE of linear elasticity with finite elements in Section \ref{FinEl}. Section \ref{lagrange} briefly recalls how to use the adjoint equation in the calculation of shape gradients following a first discretize then adjoin approach. In Section \ref{comput} we numerically calculate shape derivatives for a bended rod on the finite element mesh and we validate the calculations with the method of finite differences. We show for the given example that the geometric flows constructed from the shape derivatives, using also a suitable mesh morphing and a volume constraint, result in the optimal configuration given by the straight rod. In a second case, we study a joint connecting two levels of height, there is no optimal solution that is easily guessed. We again obtain stable geometry flows and a considerable reduction of failure probability for the numerically converged solution. In the final Section \ref{out} we give our conclusions and an outlook to future research.

\section{Failure probabilities for ceramic structures}
\label{FailProbs}

\subsection{The elasticity PDE}
\label{elasticityPDE}
Let $\Omega\subseteq\mathbb{R}^d$, $d=2,3$, be a domain with Lipschitz boundary $\partial\Omega$. It is assumed to be filled with ceramic material. 
$\Omega$ represents the ceramic component in its initial, force free state. Furthermore, we assume that the boundary $\partial\Omega$ can be divided into three different parts

\begin{align}
 \partial\Omega=\overline{\partial\Omega}_D\cup \overline{\partial\Omega}_{N_{fixed}}\cup \overline{\partial\Omega}_{N_{free}}.
\end{align}

Here, $\partial\Omega_D$ is the part of the boundary where Dirichlet-boundary conditions hold. It is supposed to be fixed. 
$\partial\Omega_{N_{fixed}}$ is the part where surface forces may act.  And finally $\partial\Omega_{N_{free}}$ is the part of the boundary 
which can be modified in order to optimally comply with the design objective 'reliability', as explained below. 
For technical reasons following from the application, the free boundary is assumed to be force-free.\\
Forces may act on the object with the shape given by $\Omega$. The volume force is represented by a function $f\in L^2(\Omega, \mathbb{R}^d)$, the surface force by a function 
$g\in L^2(\partial\Omega_N, \mathbb{R}^d)$.
In our example, the volume force $f$ represents the force of gravity, the surface force $g$ represents the tensile load.\\
Furthermore, $u\in H^{1}(\Omega, \mathbb{R}^d)$ describes the displacement caused by the acting forces. 
Here, for $m\in [1,\infty)$, $H^{1}(\Omega, \mathbb{R}^d)$ stands for the Sobolev space on $\Omega$ of once weakly differentiable $L^2(\Omega,\mathbb{R}^d)$ functions with weak derivatives in $L^2(\Omega,\mathbb{R}^{d,d})$. The linear strain tensor is given by 
$\varepsilon (u):=\frac{1}{2}\left(\nabla u+\nabla u^T\right)\in L^2(\Omega, \mathbb{R}^{d,d})$ and hence the stress tensor 
$\sigma(u)=\lambda\spur (\varepsilon(u))I+2\mu \varepsilon(u)\in L^2(\Omega, \mathbb{R}^{d,d})$, where $\lambda,\mu>0$ are the Lam\'e constants derived from Young's modulus $\texttt{E}$ and Poisson's ratio $\nu$ by $\lambda=\frac{\nu\texttt{E}}{(1+\nu)(1-2\nu)}$ and $\mu=\frac{\texttt{E}}{2(1+\nu)}$.\\
Let $H^{1}_D(\Omega,\mathbb{R}^d)$ denote the Sobolev space above with zero boundary conditions on the Dirichlet-part of the boundary. 
Under the given conditions, following from Korn's inequality on $H^{1}(\Omega;\mathbb{R}^ 3)$ \cite{duran} and a  the  Lax-Milgram theorem, the linear elasticity PDE
\begin{align}
\label{eqa:absPDE}
B(u,v)=L(v)\text{, }\forall v\in H^{1}_{D}(\Omega,\mathbb{R}^d),
\end{align}
possesses a unique weak solution $u\in H^{1}(\Omega,\mathbb{R}^d)$.  The bilinear form $B(u,v)$ and the linear functional $L(v)$ are given by
\begin{align}
\label{eqa:BandL}
\begin{split}
 B(u,v)&=\intOm[\sigma (u) : \varepsilon (v)]=\lambda\intOm[\nabla\cdot u\nabla\cdot v]+2\mu\intOm[\varepsilon (u) :\varepsilon (v)]\\
L(v) &=\int\limits_{\Omega}f\cdot v\,dx+\int\limits_{\partial\Omega_N}g\cdot v\,dA
\end{split}
\end{align}

\subsection{Survival probabilities from linear fracture mechanics}
\label{ProbFun}

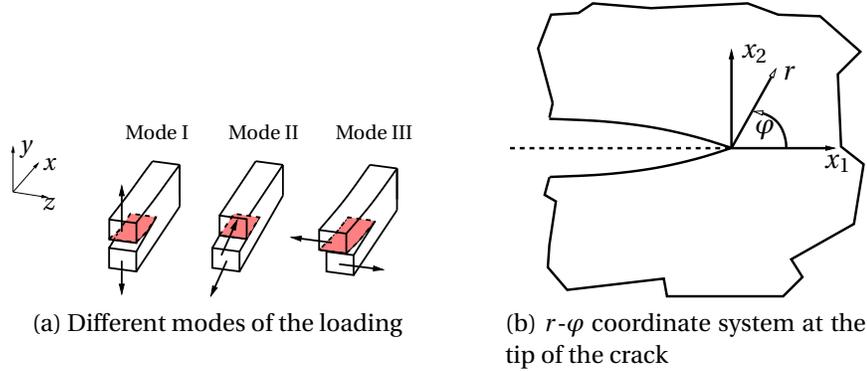
\begin{figure}[!htb]
	\begin{center}
		\subfloat[Different modes of the loading\label{fig:modes}]{\resizebox{0.4\textwidth}{!}{\input{Modus_I_II_III.pspdftex}}}
		\hspace { 1cm }
		\subfloat[$r$-$\varphi$ coordinate system at the tip of the crack\label{fig:process} ]
		{\resizebox{0.35\textwidth}{!}{\input{Fracture_mechanics_r_phi_coords_1.pspdftex}}}
	\end{center}
	\caption{Visualizations from \cite{bolten}}
\end{figure}

The full deduction of the functional describing the survival probability of a ceramic component is given in \cite{bolten}. Here we want to derive an introduction for basic understanding only.\\
To derive the objective functional we first need to study the problem in more detail. We want to optimize the reliability for a ceramic body $\Omega$ 
by minimizing its probability of failure under one given tensile load $\sigma_n$. Failure means that a fracture occurs under the tensile load. 
Therefore the question is what that probability depends on. Ceramic is produced in a process called sintering. From this process, small flaws arise in the material, 
which in the first place have no influence on the quality of the material. But under load, these flaws may become the initial point of a rupture.
To understand the behavior of this rupture, it is necessary to understand its structure and the basic principles of its growth.\\
There are three different crack opening modes. These modes are visualized in Figure \ref{fig:modes}. 
Considering the visualization of the modes, it is obvious that in the plane we only have mode I and mode II. They relate to different loadings, 
where obviously the first opening mode relates to tensile and compressive load. 

In linear fracture mechanics, the three dimensional stress field close to a crack in a two-dimensional plane close to the tip of the crack is of the form 
\begin{align}
\sigma = \frac{1}{\sqrt{2\pi r}}\{K_I\tilde{\sigma}^I(\phi )+K_{II}\tilde{\sigma}^{II}(\phi )+K_{III}\tilde{\sigma}^{III}(\phi )\}+\text{regular terms,}
\end{align}
where $r$ is the distance to the crack front and $\phi$ the angle of the shortest connection point considered to the crack front with the crack plane 
(see Figure \ref{fig:process}). Experimental evidence shows that $K_I$ is most relevant for 
the failure of ceramic structures (\cite{bruecknerfoit}) under tensile load, so we will concentrate on this.
Note that in the two-dimensional case the formula for sigma above only consists of the first two parts, as in the plane mode III load does not exist.\\
To find the functional aimed for we need to model the cracks first, as it depends on them.  We assume them to be "penny shaped" \cite{gross}. 
As a consequence, a crack can be fully described by the three properties location, orientation and radius.  
As there is no indication that one of these properties is determined by the sintering process we assume them to be arbitrary. Hence, a crack is identified 
by its configuration
\begin{align}
 (x,a,\mu )\in \left(\bar{\Omega}\times (0,\infty )\times S^{d-1}\right):=\mathpzc{C},
\end{align}
with $S^{d-1}$ the unit sphere in $\mathbb{R}^d$ and $x$ and $\mu$ are uniformly distributed on $\bar{\Omega}$ and $S^{d-1}$, respectively. The distribution of $a$ will be discussed later on.
We call $\mathpzc{C}$ the crack configuration space.\\
With this and considering the tensile load $\sigma_n$ in a normal direction of the stress plane, one obtains
\begin{align}
 K_I:=\frac{2}{\pi}\sigma_n\sqrt{\pi a}.
\end{align}
A crack becomes critical, i.e. a fracture occurs, if $K_I$ exceeds a critical value $K_{I_c}$. Obviously we can neglect compressive load, that is negative values for $\sigma_n$. With this and following \cite{bolten}, we set
\begin{align}
\sigma_n:=(n\cdot\sigma (Du)n)^+=\max \{n\cdot\sigma (Du)n,0\}.
\end{align}
Now we have determined $K_I$, we can define the set of critical configurations that lead to failure by $A_C:=A_C(\Omega,Du)=\{(x,a,\mu)\in\mathpzc{C}:K_I(a,\sigma_n(x))>K_{I_C}\}$. 
Hence we want to minimize the probability of $A_c$ containing at least one flaw.

If we assume that the distribution of cracks in different parts of the component is statistically independent, we can conclude with \cite{watanabe} and \cite{kallenberg}[Corollary 7.4.] that the random number of cracks $N(A)$ of some measurable subset of the configuration space $A\subseteq \mathpzc{C}$ is Poisson distributed, i.e. $N(A)$ is a Poisson point process (PPP). It holds that $\mathbb{P}(N(A)=n)=e^{-\nu(A)}\frac{\nu(A)^n}{n!}\sim Po(\nu(A))$, with intensity measure $\nu: \mathpzc{C}\rightarrow\mathbb{R}$. 
As mentioned before, the component fails if $N(A_C)\geq 1$.  With this, we can give the survival probability of the component $\Omega$, given in the displacement field $u\in H^1(\Omega,\mathbb{R}^d)$ as
 \begin{align}
 \label{eqa:PSurvival}
 p_s(\Omega |Du)=P(N(A_c(\Omega ,Du))=0)=\exp\{-\nu (A_c(\Omega, Du))\}.
 \end{align}
Under the assumption that cracks are statistically homogeneously distributed throughout the material and that the crack orientation is isotropic, it follows that 
 \begin{align}
 \nu=dx\otimes\nu_a\otimes \frac{dn}{2\pi^{\frac{d}{2}}/\Gamma\left(\frac{d}{2}\right)},
\end{align}  
with $dx$ the Lebesgue measure on $\mathbb{R}^d$, $dn$ the surface measure on $S^{d-1}$ and a certain positive Radon measure $\nu_a$ on $(\mathbb{R}_+,\mathcal{B}(\mathbb{R}_+))$ modelling the occurrence of cracks of length $a$. 
As mentioned before, there is a critical value $K_{I_c}$ which should not be exceeded in order to avoid failure.  Therefore, $A_c$ only contains configurations with a radius $a$ such that $K_I(a)> K_{I_c}$. This is true for all $a_c>\frac{\pi}{4}\left(\frac{K_{I_c}}{\sigma_n}\right)^2$. Due to this considerations, the intensity measure can be evaluated on the critical set as follows
\begin{align}
\nu (A_c(\Omega, Du))=\frac{\Gamma(\frac{d}{2})}{2\pi^{\frac{d}{2}}}\int\limits_{\Omega}\int\limits_{S^{d-1}}\int\limits_{a_c}^{\infty}d\nu_a(a)dndx.
\end{align}
Assuming that $d\nu (a)= c\cdot a^{-\tilde{m}}da$, with a certain constant $c>0$ and $\tilde{m}>1$, we calculate the inner integral as follows
\begin{align}
\int\limits_{a_c}^{\infty}d\nu_a(a)=\tilde{c}\left(\frac{\pi}{4}\left(\frac{K_{I_c}}{\sigma_n}\right)^2\right)^{-\tilde{m}+1}.
\end{align}
With this and setting $m:=2(\tilde{m}-1)$, with the assumption that $\tilde{m}\geq \frac{3}{2}$ holds and in assembling  all constant values in the positive constant $\left(\frac{1}{\sigma_0^m}\right)$ we find our objective functional
\begin{align}
\label{eqa:ObFun}
J(\Omega ,Du):=\nu (A_c(\Omega,Du))=\frac{\Gamma(\frac{d}{2})}{2\pi^{\frac{d}{2}}}\int\limits_{\Omega}\int\limits_{S^{d-1}}\left(\frac{\sigma_n}{\sigma_0}\right)^m dndx.
\end{align}
Let us now consider the situation, where in \eqref{eqa:BandL} we can neglect the volume force $f$ and we rescale the surface force $g$ with a constant factor $F>0$. As \eqref{eqa:absPDE} is linear, $u$, $Du$, $\sigma$ and $ \sigma_n$ all are scaled by the same factor $F$. Inserting this into \eqref{eqa:ObFun} we see that the probability of survival \eqref{eqa:PSurvival} as a function of the load parameter $F$ follows a Weibull distribution with Weibull module $m$ and scale parameter $\eta(\Omega,g)=J(\Omega,Du)^ {-\frac{1}{m}}$, where $u$ corresponds to the load scale $F=1$,
\begin{align}
p_s(F,\Omega|g)=\exp\{- J(\Omega,Du) F^ m\}=e^{- J(\Omega,Du) F^ m}=e^{-\left(\frac{F}{\eta(\Omega,g)}\right)^m}.
\end{align} 
This corresponds to the statistical strength of brittle materials as described in \cite{weibull}. Note that if we suitably normalize $g$ such that a force of one Newton is executed on the structure, $F$ can actually be interpreted as the acting force in Newton. The values for the parameter $\sigma_0$ are taken from \cite{baeker}. However we note that these values have been calibrated with a simplified model and some rescaling is needed to to fit them to our model.
\section{Discretisation via finite elements}
\label{FinEl}
To calculate the shape gradient we first need to discretise our problem with finite elements.

\subsection{Discretisation of the linear elasticity equation}
\label{FinEl-PDE}
Recall the linear elasticity PDE \eqref{eqa:absPDE} and \eqref{eqa:BandL} from section \ref{elasticityPDE}.
We discretize the PDE with standard Lagrange finite elements \cite[cf.][ch. 2.5]{braess}. For the calculation of the integrals, numerical quadrature is used. In a first step, 
$\Omega$ is partitioned by a finite mesh $\mathpzc{T}_h$ represented by the $N$ grid points $X=\{X_1,...,X_N\}$. This mesh gives as well $N_{el}$ (Lagrange) 
finite elements $\{K,\Pi (K),\Sigma (K)\}$ with $n_{sh}$ local shape functions $\theta_{K,k}\in\Pi (K)$ which are defined by the nodes $X_1^K,...,X_{n_{sh}}^K\in K$ 
with $X=\bigcup_{K\in\mathpzc{T}_h}\{X_1^K,...,X_{n_{sh}}^K\}$ and the corresponding 
Lagrange interpolation conditions
\begin{align}
\begin{array}{lr}
\varphi_j(\theta_{K,i})=\theta_{K,i}(X_j)=\delta_{ij}, & \text{for }i,j\in\{1,...,n_{sh}\}.
\end{array}
\end{align}
We assume that there exists a reference element $\{\hat{K},\hat{\Pi},\hat{\Sigma}\}$ and a bijective transformation for each element 
$K\in\mathpzc{T}_h$ $T_K:\hat{K}\rightarrow K$ such that $\hat{\Pi}=\Pi\circ T_K$, $\hat{\theta}j=\theta\circ T_K$, $j\in\{1,...,n_{sh}\}$ and 
\begin{align}
T_K=T_K(\hat{\xi},X)=\sum\limits_{j=1}^{n_{sh}}\hat{\theta}_j(\hat{\xi})X_j^K,\text{ }\hat{\xi}\in\hat{K}.
\end{align}
To numerically calculate the integral, for each $K\in\mathpzc{T}_h$ we chose $q_l^K$ quadrature points $\hat{\xi}_l^K$ on the reference element $\hat{K}$ with weights 
$\hat{\omega}_l^K$.
We then have
\begin{equation}
\label{eqa:Bdisc}
\begin{split}
B(u,v)  &=\lambda\sum\limits_{K\in\mathpzc{T}_h}\int\limits_K\nabla\cdot u\,\nabla\cdot vdx+2\mu\sum\limits_{K\in\mathpzc{T}_h}\int\limits_K\varepsilon (u):\varepsilon (v)\,dx\\
&=\lambda\sum\limits_{K\in\mathpzc{T}_h}\int\limits_K\nabla\cdot u(T_K(\hat{\xi}))\nabla\cdot v(T_K(\hat{\xi}))\text{det }(\hat{\nabla}T_K(\hat{\xi}))\, d\hat{\xi}\\
&+2\mu\sum\limits_{K\in\mathpzc{T}_h}\int\limits_K\varepsilon (u(T_K(\hat{\xi}))):\varepsilon (v(T_K(\hat{\xi})))\text{det }(\hat{\nabla}T_K(\hat{\xi}))\, d\hat{\xi}\\
&\approx\lambda\sum\limits_{K\in\mathpzc{T}_h}\sum\limits_{l=1}^{q_l^K}\hat{\omega}_l^K\text{det }(\hat{\nabla}T_K(\hat{\xi}_l))\,\nabla\cdot u(T_K(\hat{\xi}_l))\,\nabla\cdot v(T_K(\hat{\xi}_l))\\
&+2\mu\sum\limits_{K\in\mathpzc{T}_h}\sum\limits_{l=1}^{q_l^K}\hat{\omega}_l^K\text{det }(\hat{\nabla}T_K(\hat{\xi}_l))\varepsilon (u(T_K(\hat{\xi}_l))):\varepsilon (v(T_K(\hat{\xi}_l))).
\end{split}
\end{equation}
For each element $K\in\mathpzc{T}$ and $\xi\in K$ we can write $u(\xi )$ in terms of the local shape functions on the reference element: 
$u(\xi )=\sum_{m=1}^{n_{sh}}u_m\hat{\theta}_m\circ T_K^{-1}(\xi)$ and hence 
\begin{align}
\nabla u(\xi )=\sum\limits_{m=1}^{n_{sh}}u_{m}\otimes (\hat{\nabla}T_K(\hat{\xi})^T)^{-1}\hat{\nabla}\theta_m(\hat{\xi}).
\label{eq:derivU}
\end{align}
With this, we instantly get 
\begin{align}
\nabla\cdot u(x)=\sum\limits_{m=1}^{n_{sh}}\spur\left( u_{m}\otimes (\hat{\nabla}T_K(\hat{\xi})^T)^{-1}\hat{\nabla}\theta_m(\hat{\xi})\right).
\end{align}
Similar to the discretization of the bilinear form, the volume force can be discretized in the following way
\begin{align}
\label{eqa:fdic}
\int_{\Omega}f\cdot v\,dx=\sum_{K\in\mathpzc{T}_h}\sum_{l_1}^{q_l^K}\hat{\omega}_l^K\det\left(\hat{\nabla}T_K(\hat{\xi}_l)\right)f(T_K(\hat{\xi}_l) )\cdot v(T_K(\hat{\xi}_l) ).
\end{align}
The surface integral over $g$ has to be treated differently. We only consider the faces $F$ of the elements that lie on $\partial\Omega$. Let $\mathpzc{N}_h$ be the collection of these specific faces.
For each $F\in\mathpzc{N}_h$, the respective element is identified by $K=K(F)\in\mathpzc{T}_h$. One can assume that there also exists a reference face $\hat{F}$ on $\hat{K}$ such that $T_{K(F)}:\hat{F}\rightarrow F$. Additionally, for the quadrature, surface quadrature points $\hat{\xi}_l^F$ and weights $\hat{\omega}_l^F$ have to be chosen, as the face possesses one dimension less than the elements. And finally, for the transformation the square root of the Gram determinant $\sqrt{\det g_F(\hat{\xi}_l^F)}$ is required instead of the determinant of the derivative of $T_K$. It is
\begin{align}
g_F(\hat{\xi})=\hat{\nabla}^F(T_K\big |_{\hat{F}})(\hat{\xi})\left(\hat{\nabla}^F(T_K\big |_{\hat{F}})\right)^T(\hat{\xi}),
\end{align}
and thus
\begin{align}
\label{eqa:gdic}
\int_{\partial\Omega}g\cdot v\,A=\sum_{F\in\mathpzc{N}_h}\sum_{l=1}^{q_l^F}\hat{\omega}_l^F\sqrt{\det g_F(\hat{\xi}_l^F)}\,g(T_{K(F)}(\hat{\xi}_l^F))\cdot v(T_{K(F)}(\hat{\xi}_l^F)).
\end{align}
The discretized equation can be rewritten in a shorter form, in terms of the global degrees of freedom $U=(u_j)_{j\in\{1,...,N\}},u_j\in\mathbb{R}^d$ and the node coordinates $X$, where it is understood that $u_j=0$ if $X_j\in\partial\Omega_D$.
Then we have 
\begin{align}
\begin{array}{rcl}
B(X)U &=& F(X),\\
B(X)_{(j,r),(k,s)} &=& B(e_r\theta_j,e_s\theta_k),\\
F_{(j,r)} &=& \int\limits_{\Omega}f\cdot e_r\theta_jdx+\int\limits_{\partial\Omega_N}g\cdot e_r\theta_jdA;
\end{array}
\label{eq:bl}
\end{align}
with $e_r, r=1,2,3$ the standard basis on $\mathbb{R}^d$.\\

\subsection{Discretization of the objective functional}
\label{FinEl-ObFun}
For a minimal example, we first want to calculate the problem in $\mathbb{R}^2$. Therefore we consider the two dimensional objective functional
\begin{align}
 J(\Omega,u)=\int_{\Omega}\int_{S^1}\left(n\cdot\sigma\left(Du\right)n)^+\right)^mdn\,dx.
\end{align}
As $S^1$ is the unit circle in $\mathbb{R}^2$, the inner integral is
\begin{align}
I(f)=\int_{0}^{2\pi}\left(\left(\cos^2(\varphi)\sigma_{11}+2\cos(\varphi)\sin(\varphi)\sigma_{12}+\sin^2(\varphi)\sigma_{22}\right)^+\right)^md\varphi.
\end{align} 
As the function to be integrated is a periodic function over its whole period, the trapezoidal rule is the method of choice as it shows exponential convergence in this case (cf \cite{weideman}). For the present integral it yields 
\begin{align}
\begin{split}
T^{(n)}(f)&=\frac{2\pi}{n}\left(\left(\sigma_{11}^+\right)^m+\sum_{i=1}^{n-1}\left(\left(\cos^2\left(\frac{i2\pi}{n}\right)\sigma_{11}\right.\right.\right.\\
&\left.\left.\left.+2\cos\left(\frac{i2\pi}{n}\right)\sin\left(\frac{i2\pi}{n}\right)\sigma_{12}+\sin^2\left(\frac{i2\pi}{n}\right)\sigma_{22}\right)^+\right)^m\right).
\end{split}
\end{align}

By replacing $I(f)$ with $T^{(n)}(f)$, we can discretize the objective functional in the following way
\begin{align}
\label{eqa:Jdisc}
\begin{split}
J(\Omega ,u) &\approx \sum\limits_{K\in\mathpzc{T}_h}\int\limits_K T^{(n)}(f)\left(\sigma \left(x\right),\varphi\right)dx\\
&= \sum\limits_{K\in\mathpzc{T}_h}\int\limits_{\hat{K}}T^{(n)}(f)\left(\sigma \left(T_K\left(\hat{x}\right)\right),\varphi\right)\text{det }\left(\hat{\nabla}T_K\left(\hat{x}\right)\right)d\hat{x}\\
&\approx \sum\limits_{K\in\mathpzc{T}_h}\sum\limits_{l=1}^{q_l^K}\hat{\omega}_l^KT^{\left(n\right)}\left(f\right)\left(\sigma \left(T_K\left(\hat{\xi}_l^K\right)\right)\right)\text{det }\left(\hat{\nabla}T_K\left(\hat{\xi}_l^K\right)\right)
\end{split}
\end{align}

\begin{align*}
&=\sum\limits_{K\in\mathpzc{T}_h}\sum\limits_{l=1}^{q_l^K}\hat{\omega}_l^K\frac{2\pi}{n}\left(\left(\sigma\left(T_K\left(\hat{\xi}_l^K\right)\right)_{11}^+\right)^m+\sum_{i=1}^{n-1}\left(\left(\cos^2\left(\frac{i2\pi}{n}\right)\sigma\left(T_K\left(\hat{\xi}_l^K\right)\right)_{11}\right.\right.\right.\\
&\left.\left.\left.+2\cos\left(\frac{i2\pi}{n}\right)\sin\left(\frac{i2\pi}{n}\right)\sigma\left(T_K\left(\hat{\xi}_l^K\right)\right)_{12}+\sin^2\left(\frac{i2\pi}{n}\right)\sigma\left(T_K\left(\hat{\xi}_l^K\right)\right)_{22}\right)^+\right)^m\right)\\
&\cdot\text{det }\left(\hat{\nabla}T_K\left(\hat{\xi}_l^K\right)\right).
\end{align*}
In the following, we will use the finite element node set $X$ as the representative for the geometric shape $\Omega$. Likewise we use the set of global degrees of freedom $U$ to encode the (approximate) displacement field $u$. We thus write $J(X,U)$ for the discretization of $J(\Omega,u)$.
\section{Discretised shape gradients}
\label{lagrange}
After discretising the objective functional we want to calculate the shape gradient. This is
\begin{align}
 \frac{dJ(X,U(X))}{dX}=\parableit{J(X,U(X))}{X}+\parableit{J(X,U(X))}{U}\parableit{U(X)}{X}.
\end{align}
As the calculation of $\parableit{U(X)}{X}$ is very costly, we consider the corresponding Lagrange function instead. This is given by
\begin{align}
 \mathcal{L}(X,U,\Lambda):=J(X,U)-\Lambda^T(B(X)U-F(X)),
\end{align}
where $\Lambda$ is the adjoint state.\\
Calculating the derivatives of the Lagrange function with respect to all three variables yields
\begin{align}
	\begin{array}{lcr}
0\overset{!}{=}\frac{\partial \mathcal{L}(X,U,\Lambda )}{\partial \Lambda}& \Leftrightarrow & B(X)U(X)=F(X),
	\end{array}
\end{align}
which gives the state equation,
\begin{align}
	\begin{array}{lcr}
0\overset{!}{=}\frac{\partial \mathcal{L}(X,U,\Lambda )}{\partial U}=\frac{\partial J(X,U)}{\partial U}-\Lambda^TB(X) &\Leftrightarrow & B(X)\Lambda =\frac{\partial J(X,U)}{\partial U},
	\end{array}
\end{align}
which is the adjoint equation. Hence, the following set of equations
\begin{align}
\begin{split} 
\frac{dJ(X,U(X))}{dX}&=\frac{\partial J(X,U)}{\partial X}+\Lambda^T\left[\frac{\partial F(X)}{\partial X}-\frac{\partial B(X)}{\partial X}U\right] \\
B^T(X)\Lambda &= \frac{\partial J(X,U)}{\partial U}\\
B(X)U(X) &= F(X)
\label{eqa:adjointEq}
\end{split}
\end{align}
gives the discretized shape gradient. 

\section{Computation and validation of shape gradients and shape flows}
\label{comput}

\subsection{Implementation}
The calculation of the shape gradient using the adjoint formalism \eqref{eqa:adjointEq} requires the numerical calculation of the state $U$, $\frac{\partial J(X,U)}{\partial U}$, of the adjoint state $\Lambda$ and the calculation of $\frac{\partial J(X,U)}{\partial X}$, $\frac{\partial B(X)}{\partial X}$ and $\frac{\partial F(X)}{\partial X}$. This can be done by somewhat lengthy but straight forward calculations based on \eqref{eqa:Jdisc}, \eqref{eqa:Bdisc}, \eqref{eqa:fdic} and \ref{eqa:gdic}. 

For the implementation, all these partial derivatives are calculated locally for each local node set of the finite elements and are assembled to global objects thereafter. Note that the contractions with the adjoint state $\Lambda$ and the state $U$ have to be performed during local calculations prior to the assembly in order to keep memory requirements for the storage of $\frac{\partial F(X)}{\partial X}$ and especially $\frac{\partial B(X)}{\partial X}$ reasonably low. 

\subsection{Validation with finite differences}
\begin{figure}
	\begin{center}
		{\includegraphics[width=0.5\textwidth]{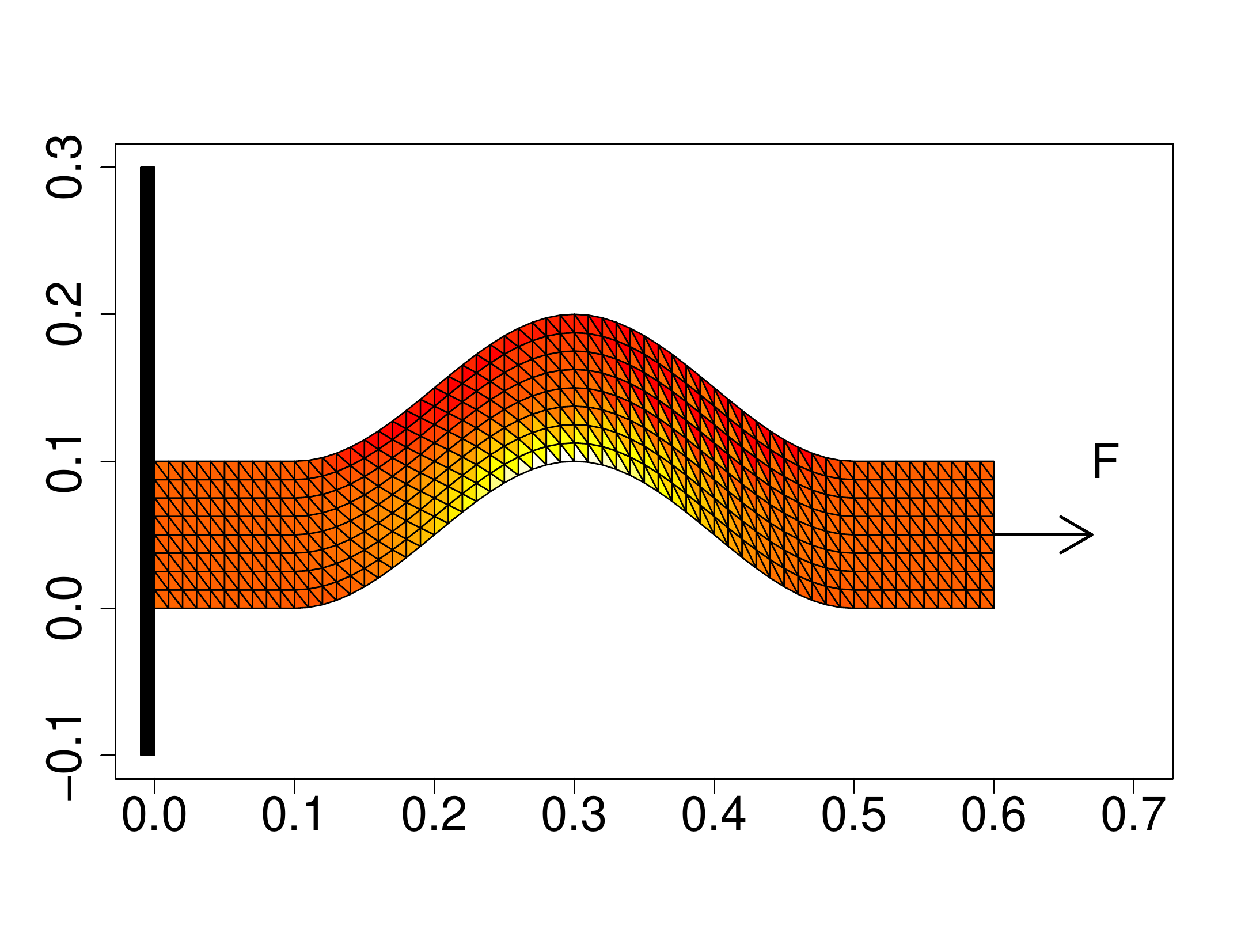}}
	\end{center}
	\caption{ Visualization of the objective functional}
	\label{fig:ob_fun_val}
\end{figure}

We consider a simple example in $d=2$: For this purpose we generate a simple two-dimensional test object whose behavior during the optimization process is well understood. 
To work with reasonable values we set the parameters $\texttt{E}$ and $\nu$  to those of Aluminum oxide ($Al_2O_3$) ceramics. The elastic material properties can be found in \cite{aluminum}. The Weibull modulus $m$ is measured in tensile tests \cite{baeker,munz}.  As it generally depends on the technical details of the sintering process, $m$ is not a materials property that is determined by the chemistry. Technical ceramics usually comes with a Weibull module between 5 for low quality and 20 for a very controlled process. Here we choose $m=10$, which is a reasonable value and sill leads to  tractable numerics.\\
As a test object we use a rod of length 0.6m and height 0.1m. The test object is fixed on the left boundary, that is where Dirichlet-boundary conditions hold. The part of the surface where 
surface forces may act is the right boundary, that is in our model the nodes on the right edge. In our example, we suppose this boundary part to be fixed as well. It is represented by a 9x61 grid, that is divided in triangles. The rod is deformed in the middle part, see Figure \ref{fig:ob_fun_val}. 

As element type we choose linear, Lagrangian triangle elements with three local degrees of freedom located at the vertices of the triangles. 
For the interpolation of the volume force and the surface force we use two-dimensional 7-point Gauss-quadrature and one-dimensional 3-point Gauss-quadrature, respectively. 
To construct the bilinear form we use 1-point Gauss-quadrature.

To fit our purpose, we use a direct finite element solver written in \texttt{R} for the underlying elasticity equation. With the results of the solver, one can calculate the values of the objective functional on the test object \eqref{eqa:Jdisc}. One can see in the visualization in Figure \ref{fig:ob_fun_val} 
that the local intensity for the occurrence of critical cracks, i.e. the density with respect to $dx$ of $\nu_a(dx\times S^1\times [a_c\infty])$, takes the highest values in the critical area in the inner bow of the deformation, as practitioners would of course expect.

\subsection{Calculation and validation of the local and the global derivatives}
\label{vallocal}

 \begin{figure}
 	\begin{center}
 		{\includegraphics[width=0.5\textwidth]{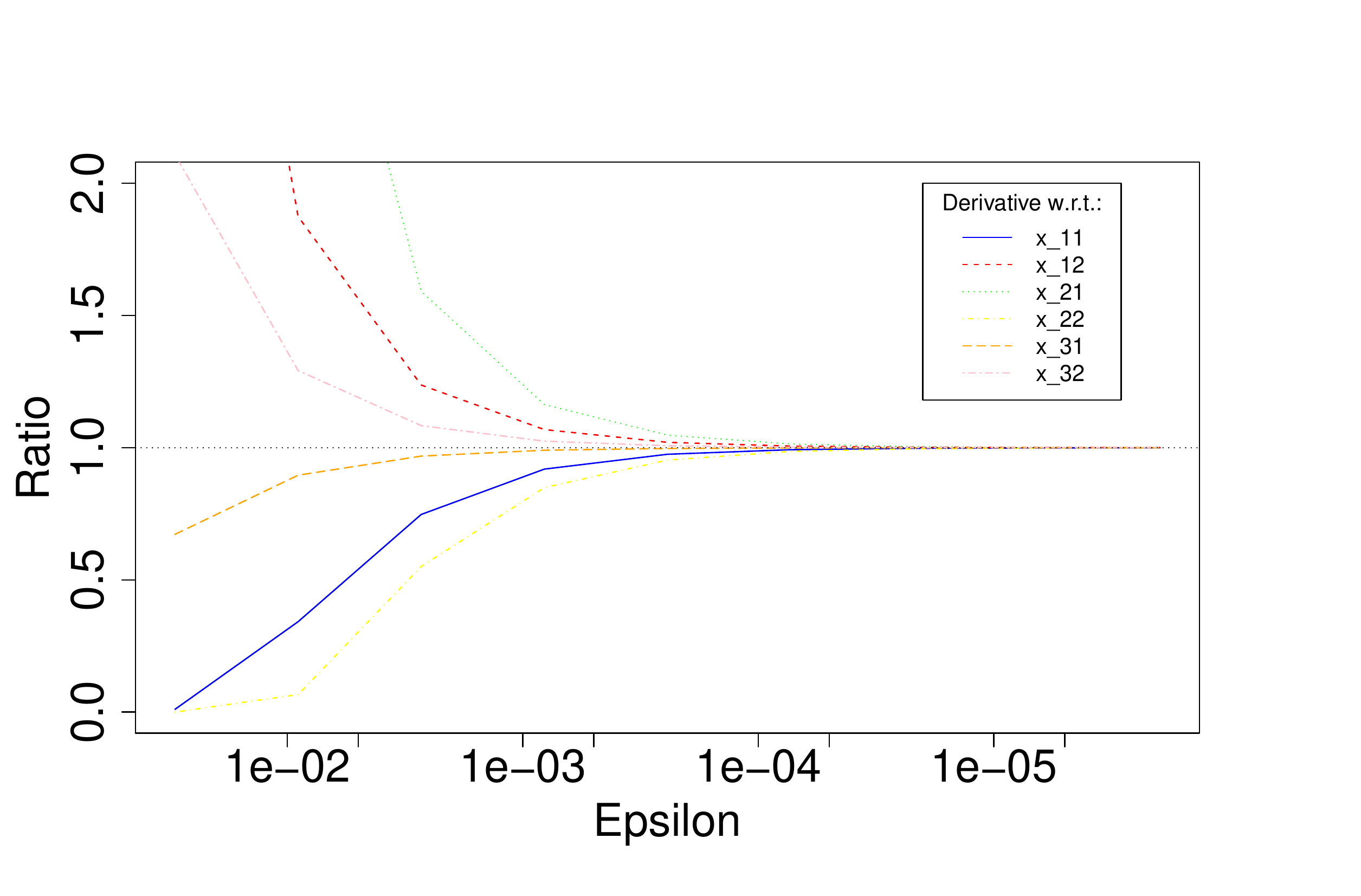}}
 	\end{center}
 	\caption{ Convergence test for the local partial derivative $\parableit{J}{X}$}
 	\label{fig:conv_dx}
 \end{figure}

\begin{figure}
	\begin{center}
	{\includegraphics[width=0.5\textwidth]{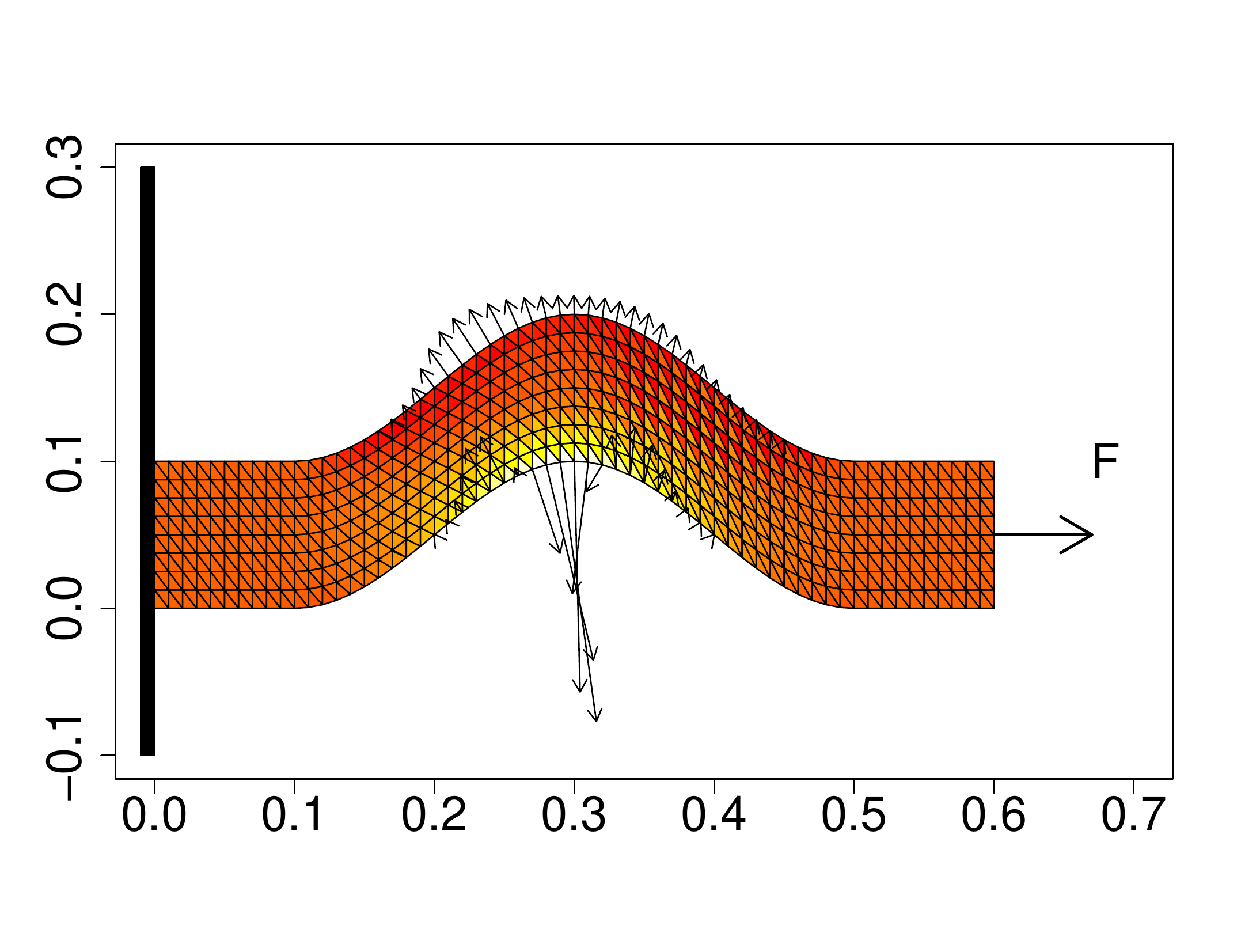}}
	\end{center}
	\caption{Visualization of the gradient $-\frac{dJ}{dx}$}
	\label{fig:finit_Diff}
\end{figure}

\begin{figure}
	\centering
	\begin{minipage}{.5\linewidth}
		\centering
		{\includegraphics[width=\linewidth]{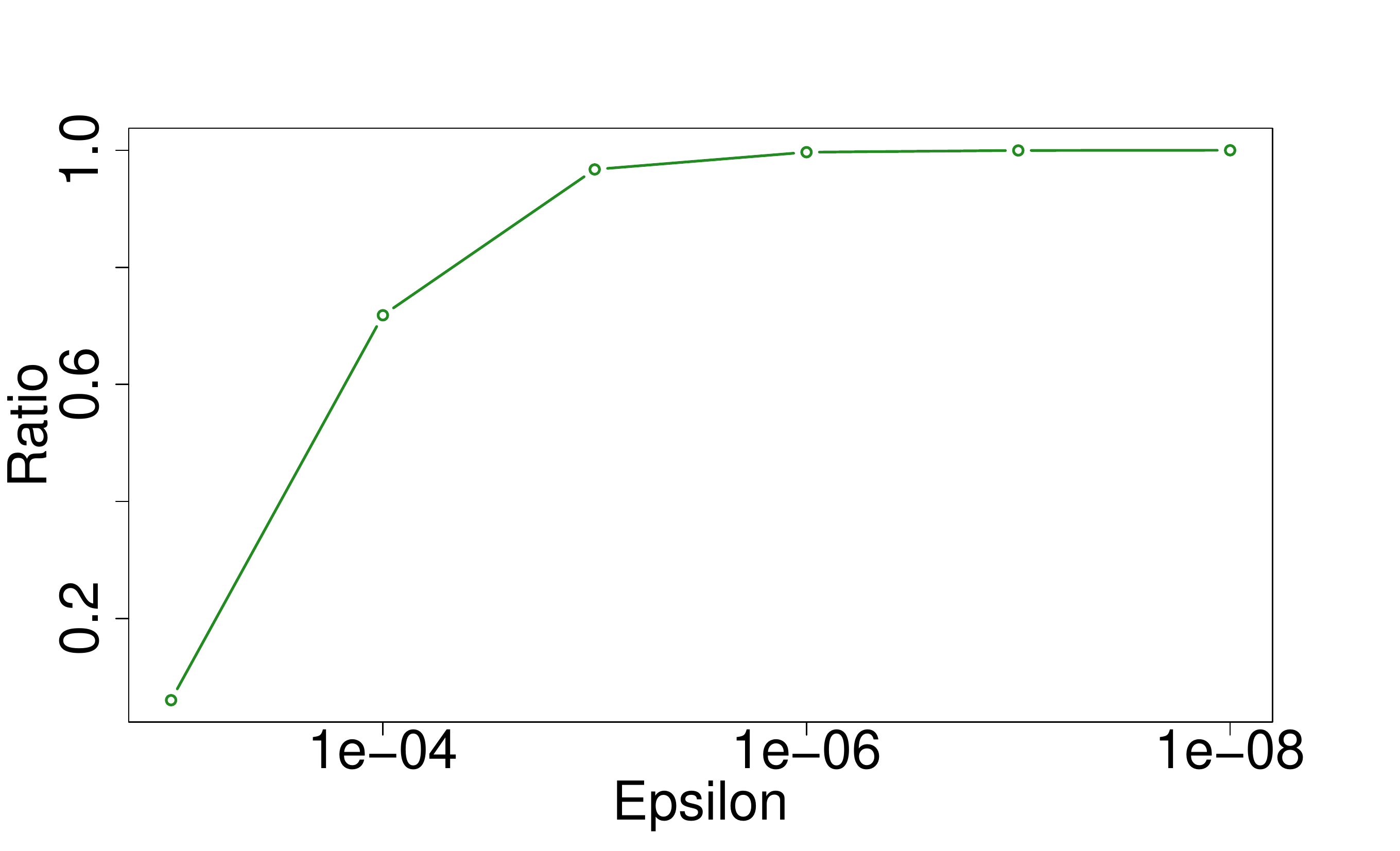}}
		\par{\vspace{0pt}}
	\end{minipage}\hfill
	\begin{minipage}{.4\linewidth}
		\centering
		\vspace{0pt}
		\begin{tabular}{lc}
			\toprule
			Epsilon & Ratio \\
			\midrule
			$10^{-3}$ & 0.0603255\\
			$10^{-4}$ & 0.7180954\\
			$10^{-5}$ & 0.9672460\\
			$10^{-6}$ & 0.9966801\\
			$10^{-7}$ & 0.9996691\\
			$10^{-8}$ & 0.9999882\\
			\bottomrule
		\end{tabular}
		\par{\vspace{0pt}}
	\end{minipage}
	\caption{Convergence of $\,||\frac{dJ(X,U)}{dX}||_2/||(\frac{J(X+\varepsilon V,U(X+\varepsilon V))-J(X,U)}{\varepsilon V})||_2$ for the two-dimensional objective functional}
	\label{fig:conv_DJDX}
\end{figure}

In Figure \ref{fig:conv_dx} the comparison of the differential quotient with the derivative of $\parableit{J}{X}$ is shown for the local derivatives on one element $K$ for five random directions.  Convergence for small differences is excellent.  For the partial derivatives, the same findings applied.    
As the node set of the test object is comparatively small, we can approximate the global shape derivative of the objective functional $\frac{d J}{d X}$ by its differential quotient for validation purposes. We calculate $\frac{J(X+\varepsilon V,U(X+\varepsilon V))-J(X,U(X))}{\varepsilon V}$, where $V$ is a random direction, with decreasing $\varepsilon$ and tested our resulting gradient against it. The results are visualized in Figure \ref{fig:conv_DJDX}. The actual gradient is visualized in Figure \ref{fig:finit_Diff}.

With this result, evidence is given that the shape gradients of the objective functional work in a reasonable way. For the mechanical interpretation, it is visible that the negative gradient tends to increase the volume of the body, which is natural, as the force is constant and so a higher volume would diminish the force acting on each node - here we have set the gravitational force $f$ to zero. The problem that one can not expect an optimal solution for the given problem, if the gravitational force $f$ is absent, will be handled via a volume constraint in the optimization process.\\ 
Apart from the visualization, we also compare the euclidean norm of the finite differences solution with decreasing epsilon with the norm of the solution of the program. This is given in Figure \ref{fig:conv_DJDX}. One can see that a very good convergence is obtained.\\
\subsection{Shape flows towards higher reliability}
\label{optim}
To obtain a first optimization process, we use geometric mesh morphing and small step sizes. I.e. we apply a mapping $\theta\to X(\theta)$ such that the interior nodes of the node set $X$ follow the movement of the surface nodes such that a reasonable mesh quality is preserved. In our case, the design parameter set $\theta$ consists of all $y$-coordinates of all nodes between the two faces on the rod in $x$-direction, on which boundary conditions are applied. Thus, the shape is (up to discretization) not restricted by this parametrization. The internal nodes with a fixed $x$ coordinate are distributed with equal distances between the respective surface nodes. The shape gradient w.r.t. the $y$-position of the surface nodes, $\frac{dJ}{d\theta}=\frac{dJ}{dX}\frac{dX}{d\theta}$, can be easily calculated as $\frac{dX}{d\theta}$ is available analytically. The optimization problem in the design parameters $\theta$ thus is 118-dimensional (there are 61 rows with two surface nodes each and two of these rows are fixed). 

\begin{figure}[!htb]
	\begin{center}
		\subfloat[\label{fig6:Bild1} ]{\includegraphics[width=0.3\textwidth]{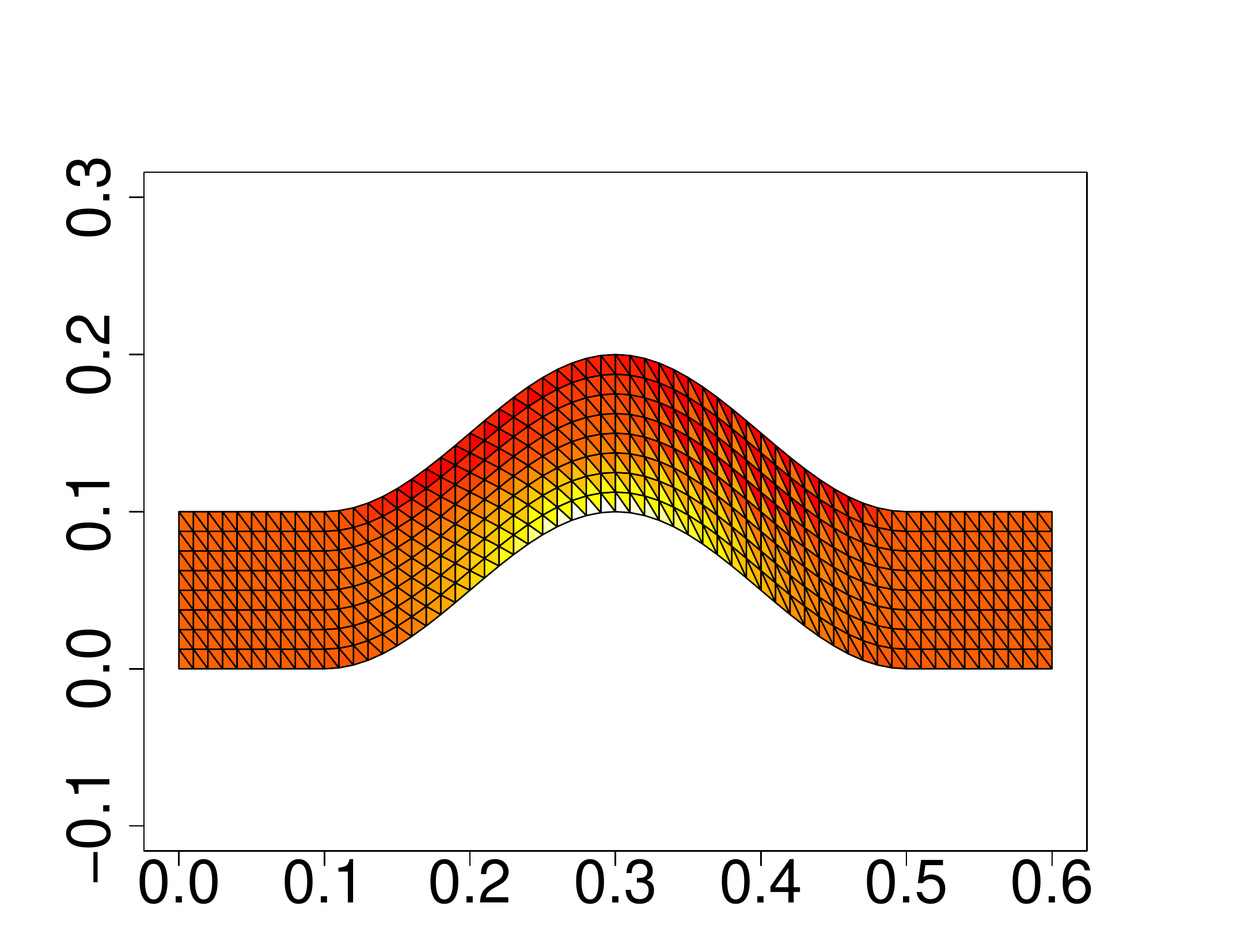}}
		\hspace { 1cm }
		\subfloat[\label{fig6:Bild2} ]{\includegraphics[width=0.3\textwidth]{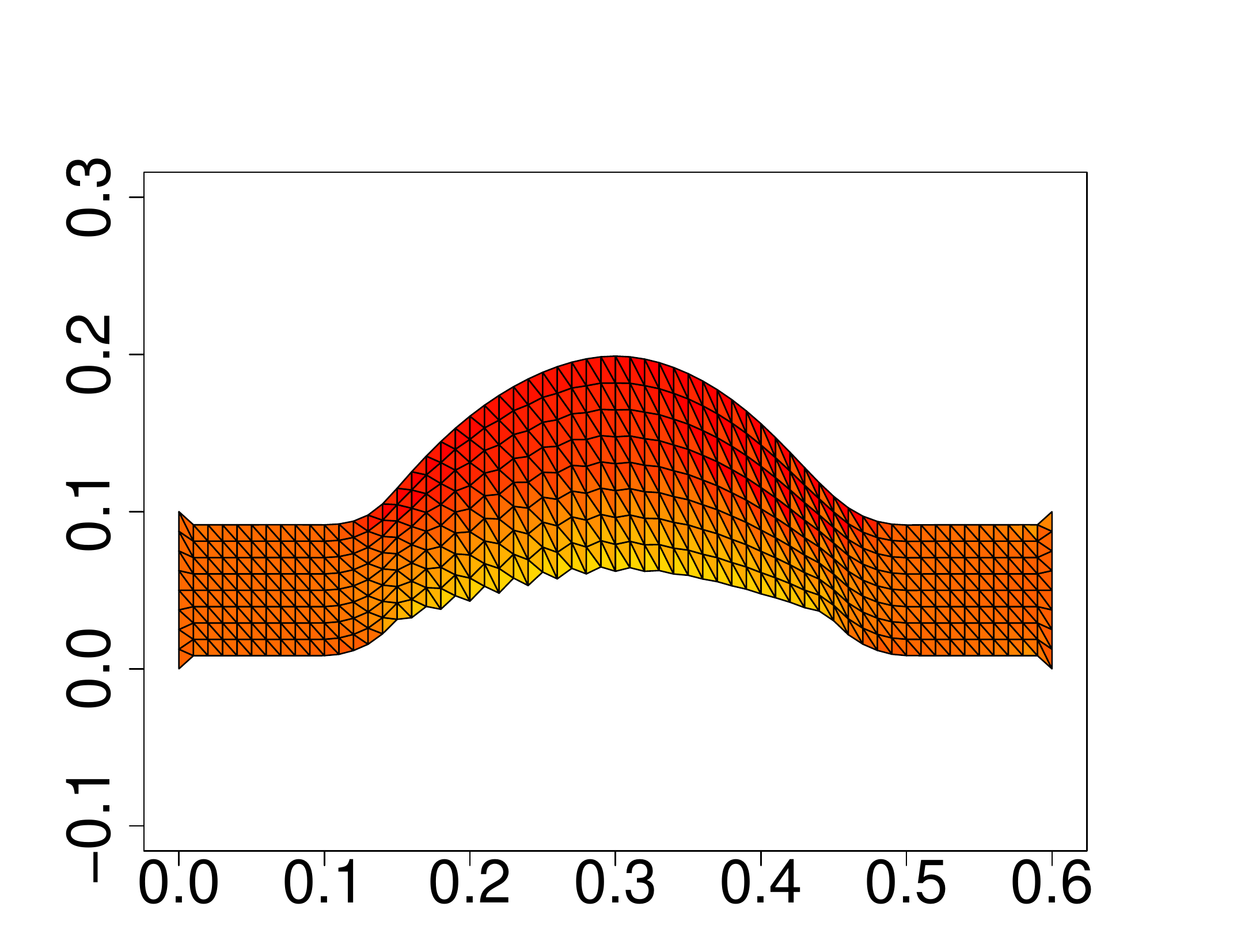}}
		\hspace { 1cm }
		\subfloat[\label{fig6:Bild3} ]{\includegraphics[width=0.3\textwidth]{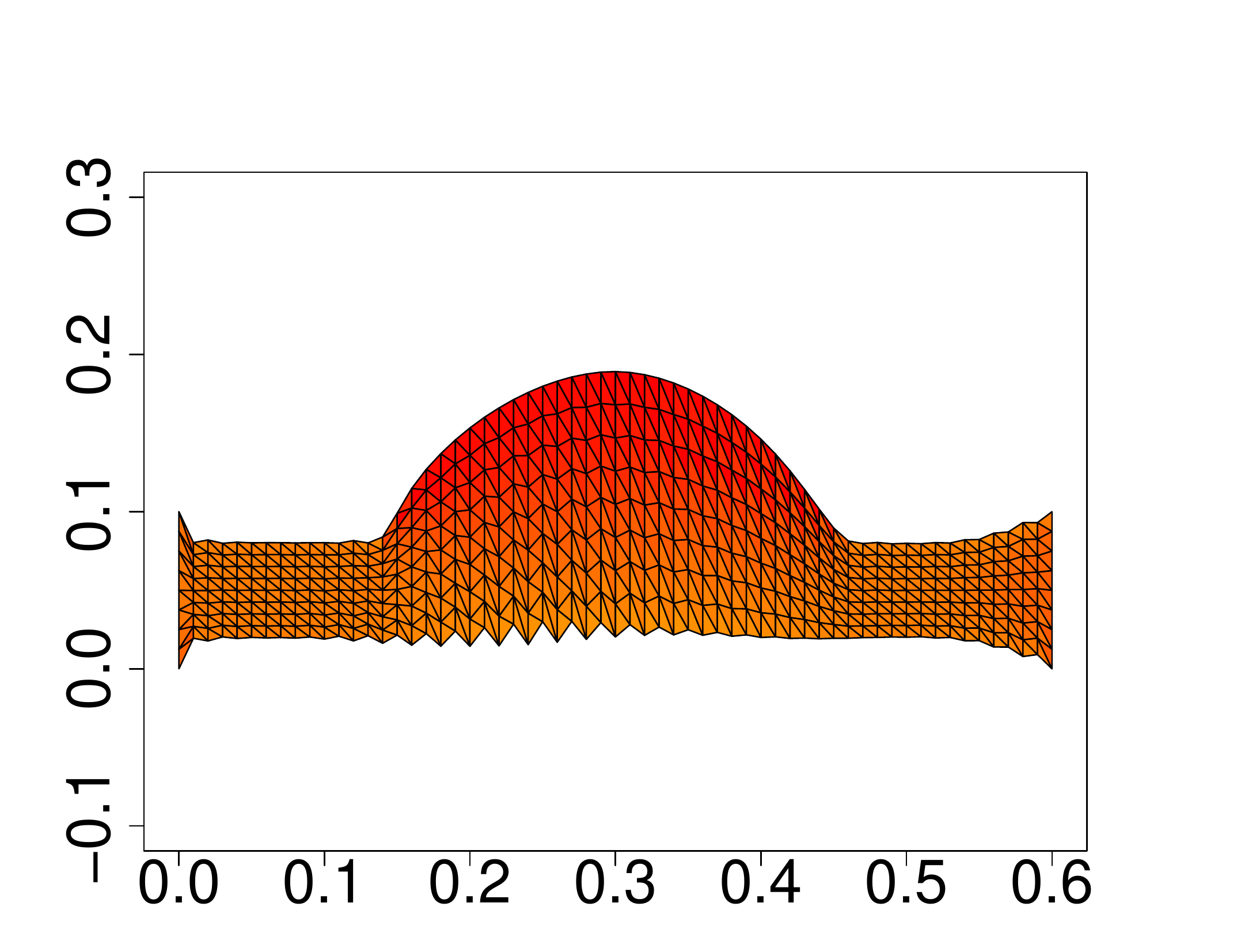}}
		\hspace { 1cm }
		\subfloat[\label{fig6:Bild4} ]{\includegraphics[width=0.3\textwidth]{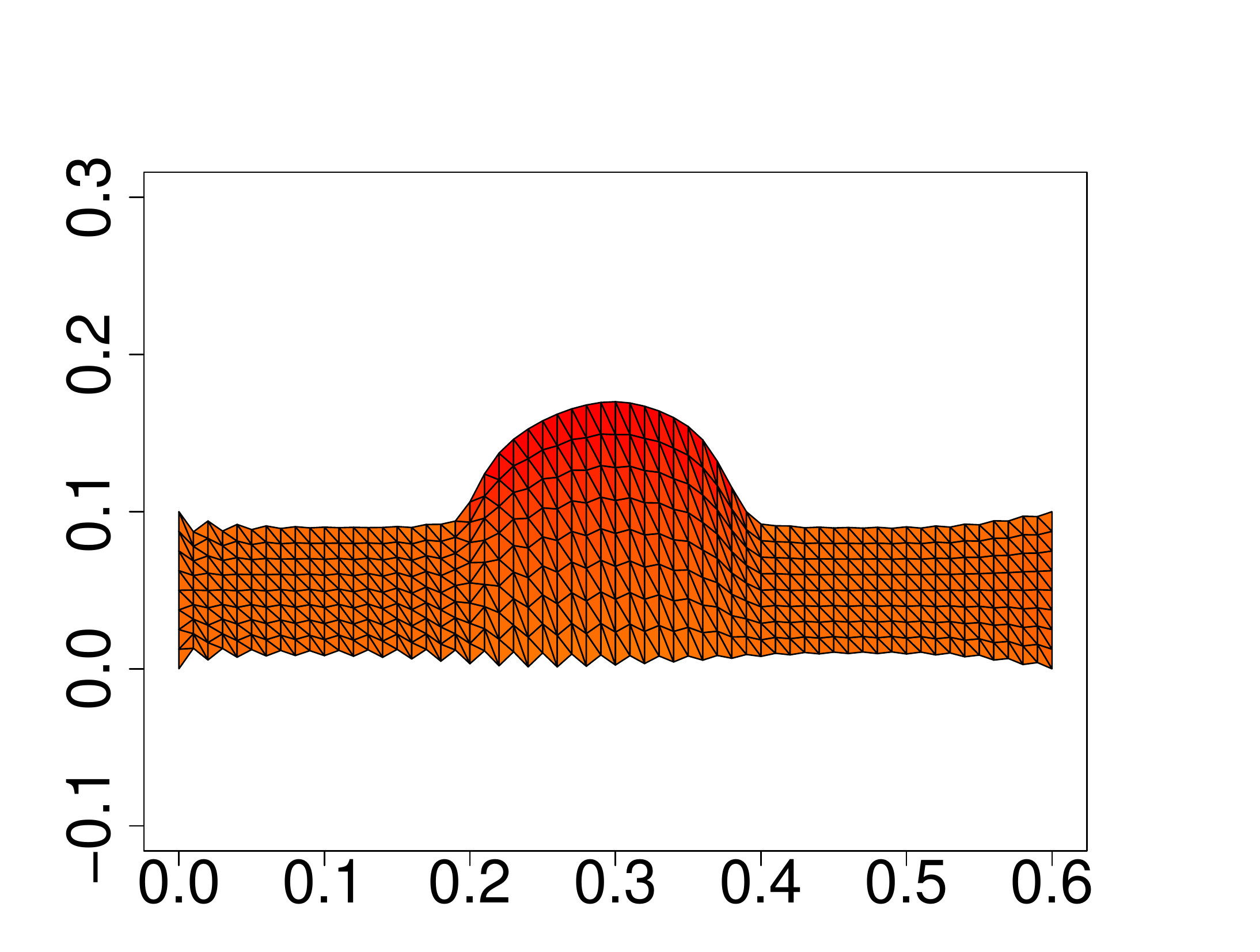}}
		\hspace { 1cm }
		\subfloat[\label{fig6:Bild5} ]{\includegraphics[width=0.3\textwidth]{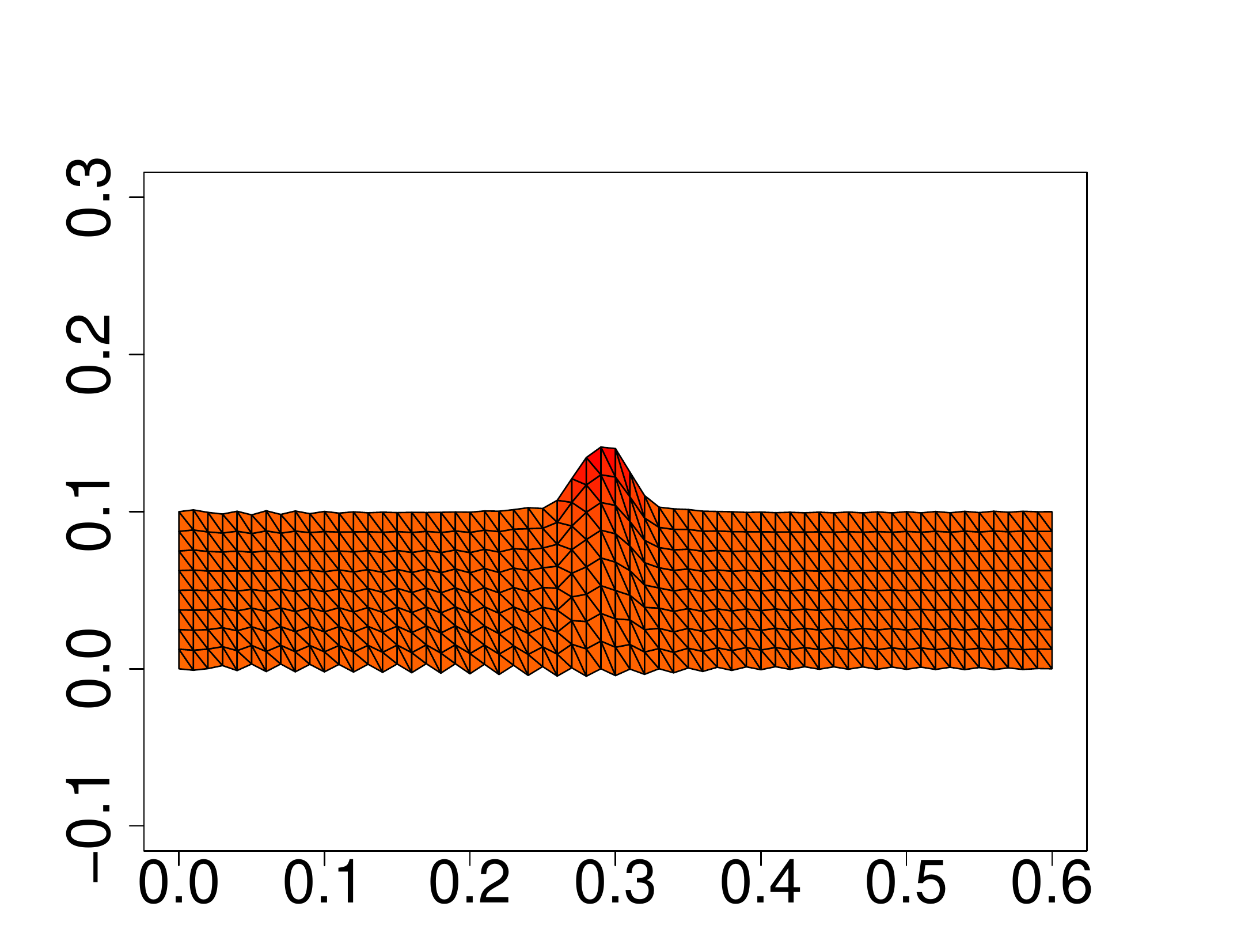}}
		\hspace { 1cm }
		\subfloat[\label{fig6:Bild6} ]{\includegraphics[width=0.3\textwidth]{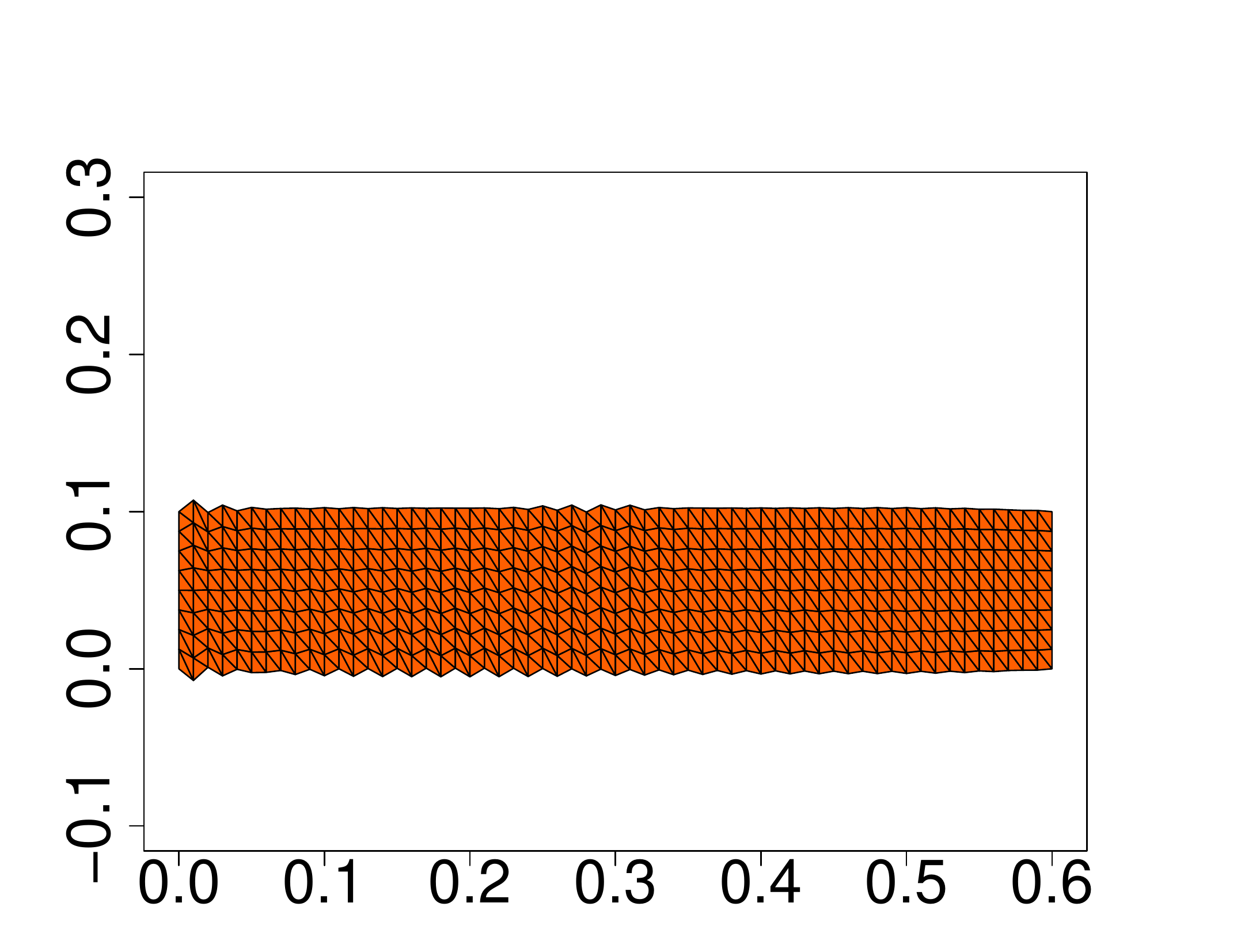}}
	\end{center}
	\caption{ The optimization procedure }
	\label{fig:optim}
\end{figure}

\begin{figure}
	\begin{center}
		{\includegraphics[width=0.5\textwidth]{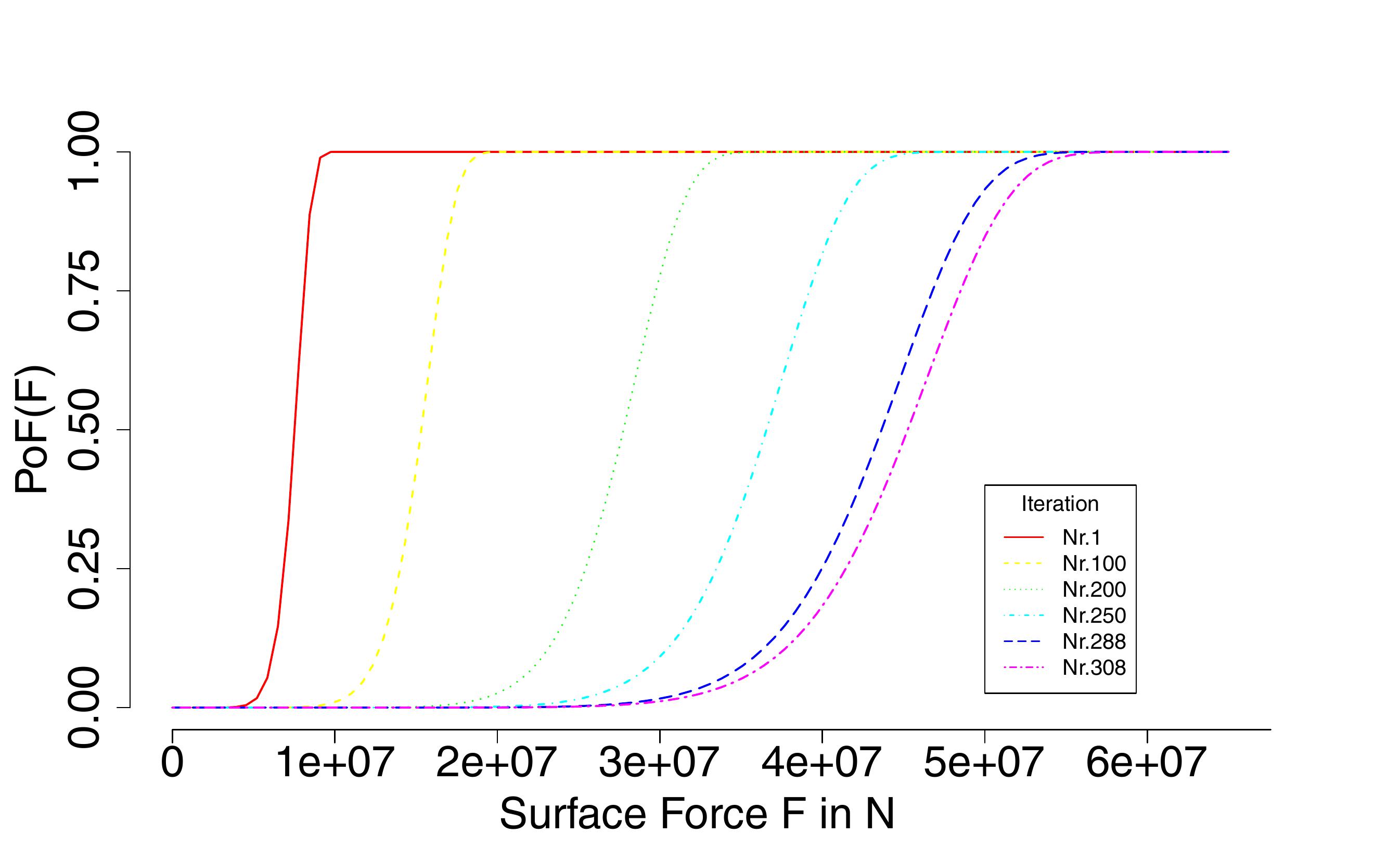}}
	\end{center}
	\caption{ Distribution of the bended rod's UTS during the optimization process}
	\label{fig:prob}
\end{figure}

To keep the volume of the body constant, we calculate the volume gradient and project it from the shape gradient. To determine the new mesh in each iteration, we calculate the new boundary nodes by $\theta_{new}=\theta_{old}-\alpha \left(\frac{dJ}{d\theta}-\frac{\partial \text{Vol}}{\partial \theta}\right)$.\\
The outcome is a procedure, converging against the straightened rod, which is visualized in Figure \ref{fig:optim}. We can also investigate the (discretised) failure probability of the component for each iteration over the load parameter $1-p_s(t,\Omega|g)$, cf. \eqref{eqa:PSurvival}, where the reference $g$ for $t=1$ is chosen such that the resulting force is one N. In Figure \ref{fig:prob} we see the distributions for some iteration steps. It is obvious that we actually decrease the probability of failure with the present procedure, which finally converges to the evident optimum.

\begin{figure}
	\begin{center}
		{\includegraphics[width=0.5\textwidth]{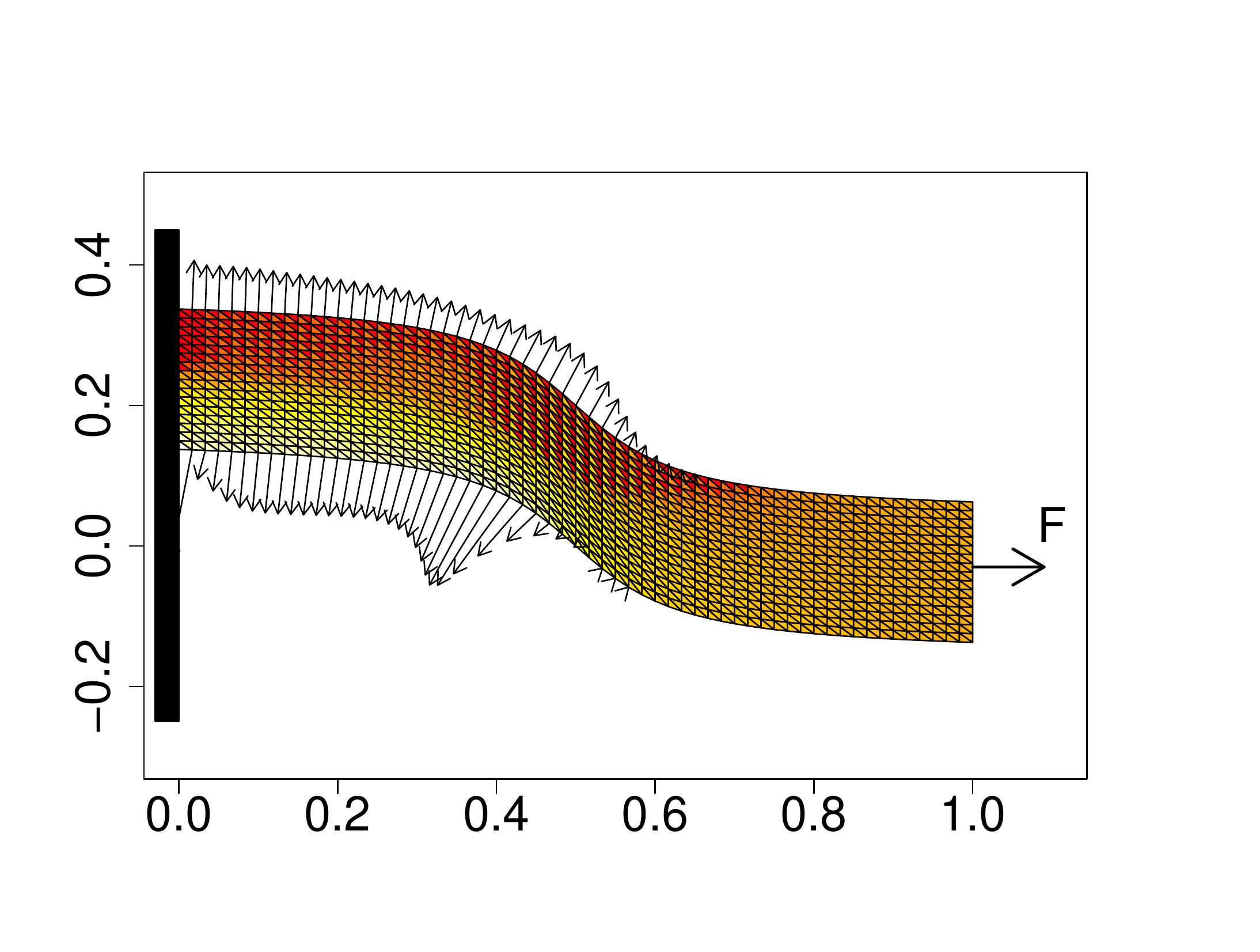}}
	\end{center}
	\caption{ Visualization of the objective functional and the gradient of the S-shaped joint}
	\label{fig:Ob_Fun_S}
\end{figure}

We also test the procedure to improve the reliability of a S-shaped joint visualized in Figure \ref{fig:Ob_Fun_S}, combined with the visualization of the gradient. It is represented by a 61x17 mesh.
For this joint it is not obvious in which direction the shape is supposed to tend to reduce the failure probability. 

\begin{figure}
	\begin{center}
		{\includegraphics[width=0.5\textwidth]{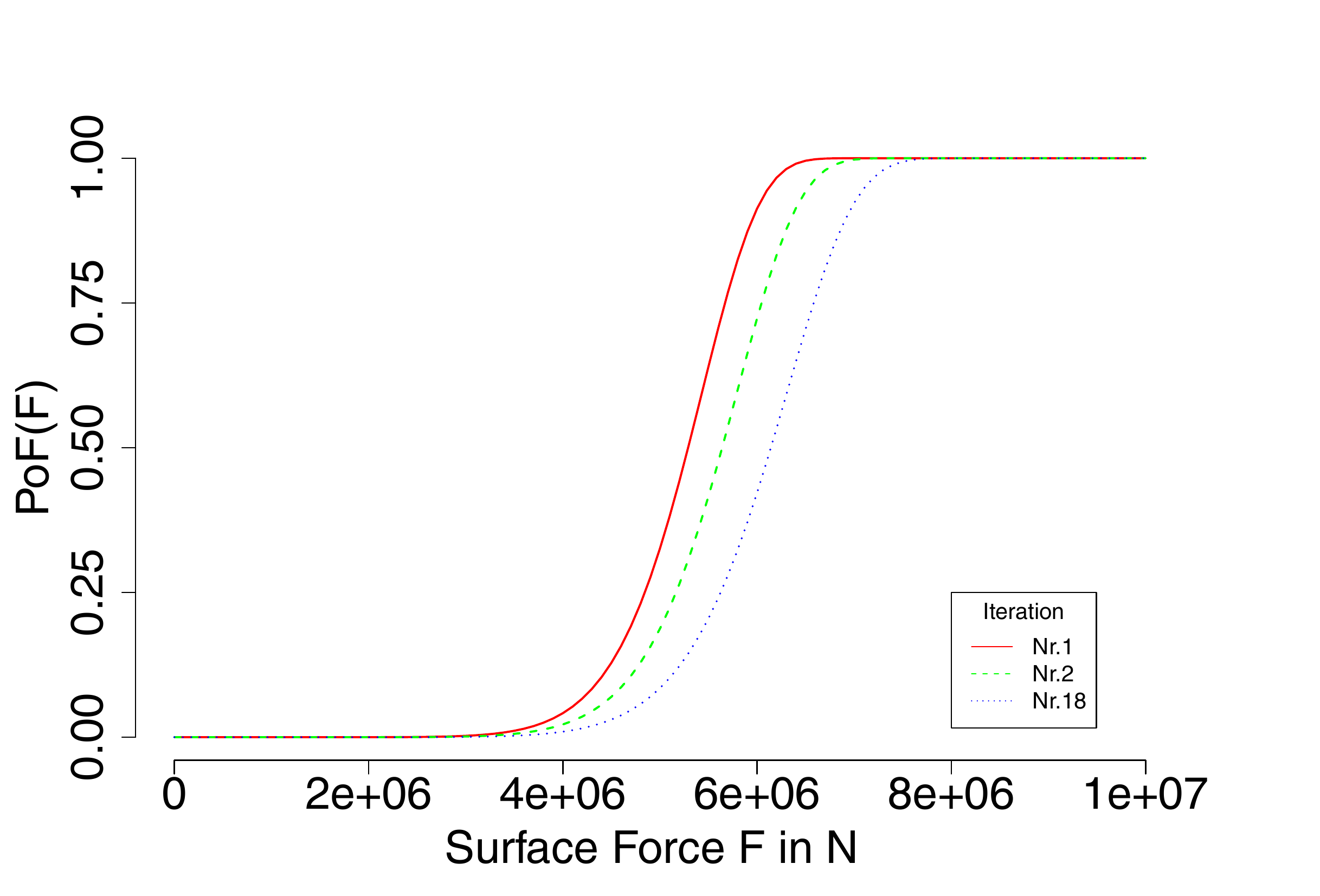}}
	\end{center}
	\caption{ Distribution the UTS of the joint for different iterations in the optimization process}
	\label{fig:prob_S}
\end{figure}

Here as well, we did apply the optimization procedure we already applied to the first example. In Figure \ref{fig:prob_S} it is visible that, even though it is not a strict descent procedure, we are able to reduce the failure probability of the component significantly.

 \begin{figure}
	\centering
	\begin{minipage}{.4\linewidth}
		\centering
		{\includegraphics[width=\linewidth]{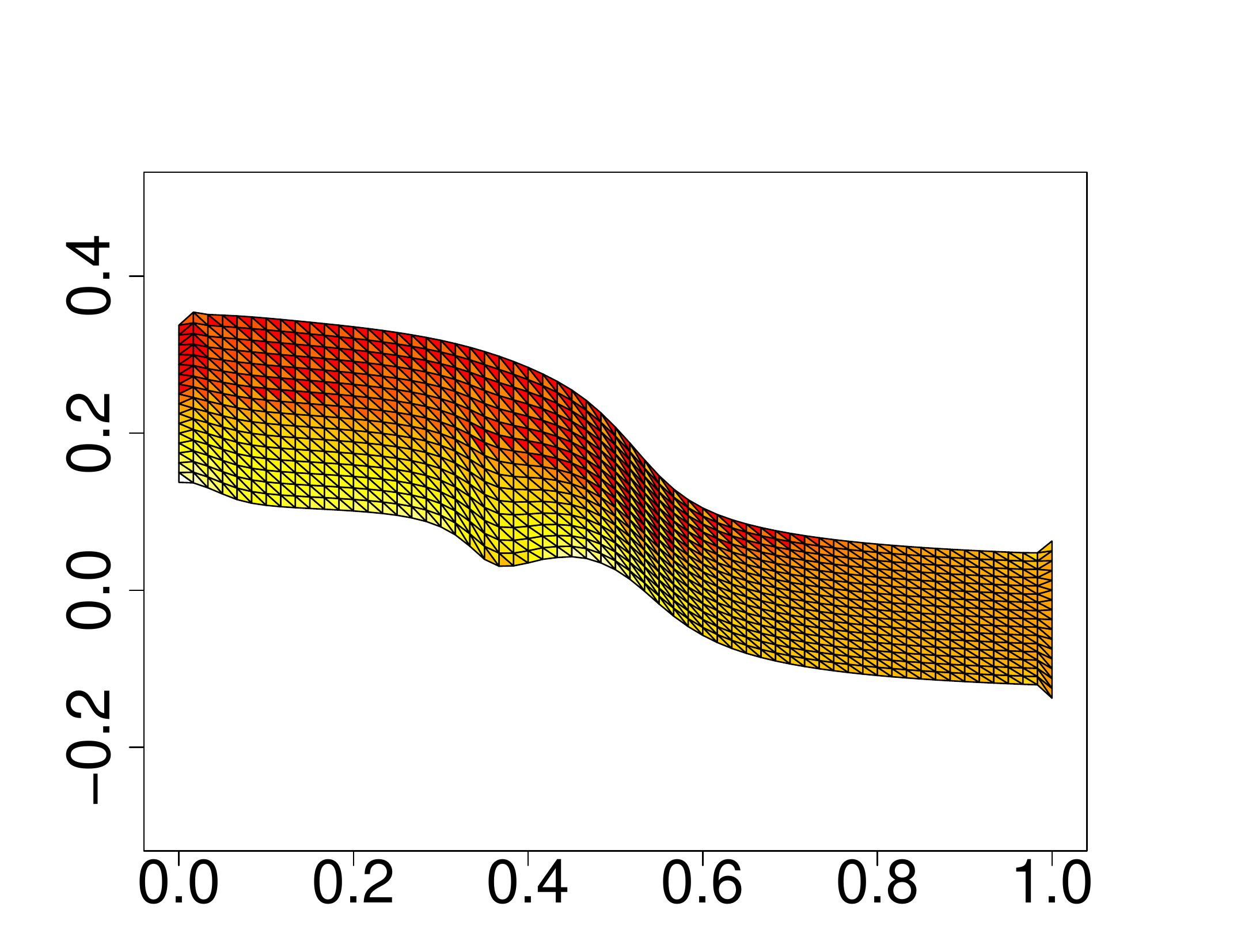}}
		\par{\vspace{0pt}}
		\caption{Iteration 2}
	\end{minipage}\hfill
	\begin{minipage}{.4\linewidth}
			\centering
		{\includegraphics[width=\linewidth]{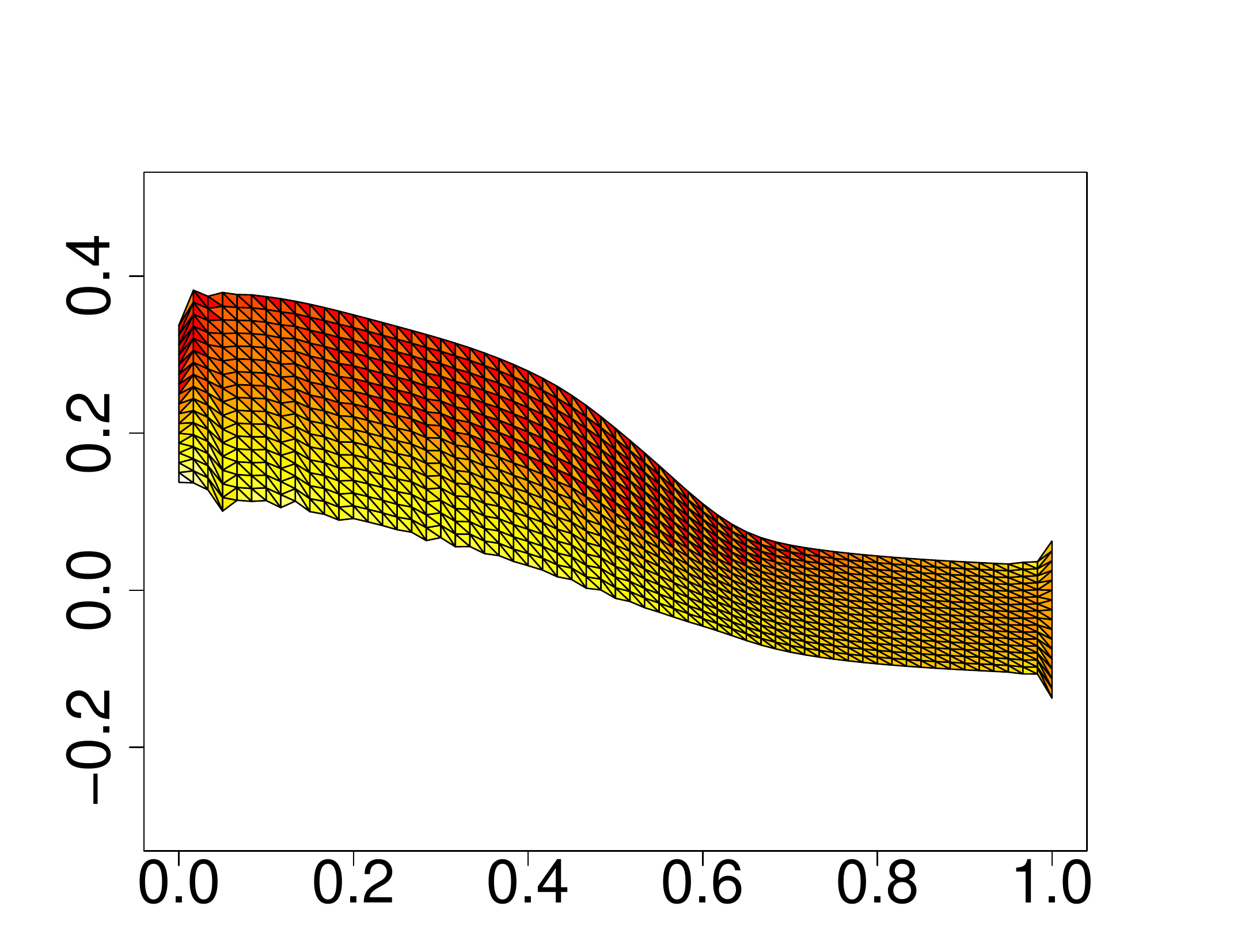}}
		\par{\vspace{0pt}}
		\caption{Iteration 18}
	\end{minipage}
	\caption{Results of the optimization of the joint}
	\label{fig:op_S}
\end{figure}

\section{Conclusion and Outlook}
\label{out}
This paper for the first time provides a numerical shape optimization algorithm that is directly connected to the  design criterion of (probabilistic) mechanical integrity. Based on a model for the probabilistic UTS \cite{batendorf,bruecknerfoit,weibull}, we computed the discretized shape gradient for the probability of failure of a ceramic component and embedded this in gradient based optimization procedure. We validated the implementation and showed that gradient flows with volume constraint are stable over major changes of the geometry. In a situation, where the optimal shape is intuitively clear, we observed convergence of the algorithm.  

While the algorithm gives quite satisfactory solutions from an intuitive engineering standpoint, there are several rather obvious research necessities from a more theoretical prospective: The first question is, in what sense the solution of the discretized problem is close to the solution of the shape optimization in not necessarily polygonal shapes. While the general strategy for such a proof is clear \cite{haslinger}, in the given context are some subtle points to keep track of: The discretized optimal shapes in our case are polygons. Given the problem of stress divergence $\sim r^{\lambda(\gamma)}$ as a function of the opening angle $\gamma$ of a reentrant corner \cite{nazarov}, the non discretized solutions on the discrete domains $\Omega_h$ will have diverging objective functionals $J(\Omega_h,u)$ if $m\lambda(\gamma)\geq d$. Given the high values of technical Weibull modules $m\approx 10\ldots 20$, this imposes a kind of discrete curvature restriction on the domains $\Omega_h$, where an error estimate is feasible. While one would expect that this problem will go away by itself if the mesh size $h$ is sufficiently fine, as the optimal configuration will only pick shapes $\Omega_h$ which stay away from the critical bound, it is somewhat tricky to cast this to a mathematical proof.  We intend to revisit this problem in the future.

The second obvious challenge  is a first adjoin and then discretize solution for the given problem. This strategy however requires the shape differentiation of an objective functional, which is not even continuous in $H^1(\Omega,\mathbb{R}^d)$. This however can be achieved using smoother shapes and elliptic regularity theory \cite{agmon,agmon2,ciarlet,gottschalk}. An adequate discretization of the continuous shape derivatives thus also requires a careful handling of surface smoothness. 

\vspace{.3cm}

\noindent\textbf{Acknowledgements:} \\
Hanno Gottschalk thanks Sebastian Schmitz form Siemens Energy, Gas Turbine Engineering Department, for interesting discussions. 


{\small
\lineskip=.1cm
\bibliographystyle{plain}
\bibliography{literatur}
}

\end{document}

%% file: Modus_I_II_III.pspdftex
\begin{picture}(0,0)%
\includegraphics{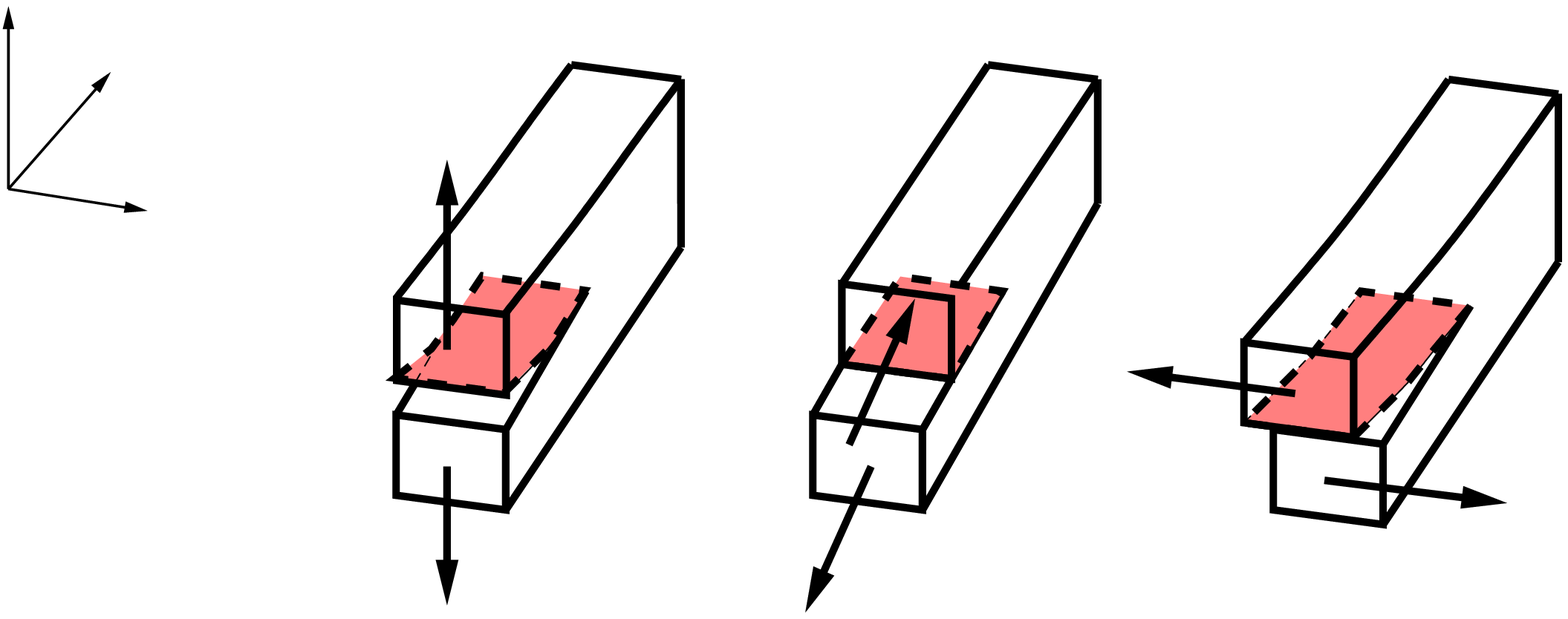}%
\end{picture}%
\setlength{\unitlength}{4144sp}%
\begingroup\makeatletter\ifx\SetFigFont\undefined%
\gdef\SetFigFont#1#2#3#4#5{%
  \reset@font\fontsize{#1}{#2pt}%
  \fontfamily{#3}\fontseries{#4}\fontshape{#5}%
  \selectfont}%
\fi\endgroup%
\begin{picture}(10246,4394)(309,-4740)
\put(1081,-1546){\makebox(0,0)[lb]{\smash{{\SetFigFont{41}{49.2}{\familydefault}{\mddefault}{\updefault}{\color[rgb]{0,0,0}$x$}%
}}}}
\put(5716,-781){\makebox(0,0)[lb]{\smash{{\SetFigFont{34}{40.8}{\familydefault}{\mddefault}{\updefault}{\color[rgb]{0,0,0}Mode II}%
}}}}
\put(8371,-781){\makebox(0,0)[lb]{\smash{{\SetFigFont{34}{40.8}{\familydefault}{\mddefault}{\updefault}{\color[rgb]{0,0,0}Mode III}%
}}}}
\put(3151,-781){\makebox(0,0)[lb]{\smash{{\SetFigFont{34}{40.8}{\familydefault}{\mddefault}{\updefault}{\color[rgb]{0,0,0}Mode I}%
}}}}
\put(1081,-2536){\makebox(0,0)[lb]{\smash{{\SetFigFont{41}{49.2}{\familydefault}{\mddefault}{\updefault}{\color[rgb]{0,0,0}$z$}%
}}}}
\put(496,-1096){\makebox(0,0)[lb]{\smash{{\SetFigFont{41}{49.2}{\familydefault}{\mddefault}{\updefault}{\color[rgb]{0,0,0}$y$}%
}}}}
\end{picture}%

%% file: Fracture_mechanics_r_phi_coords_1.pspdftex
\begin{picture}(0,0)%
\includegraphics{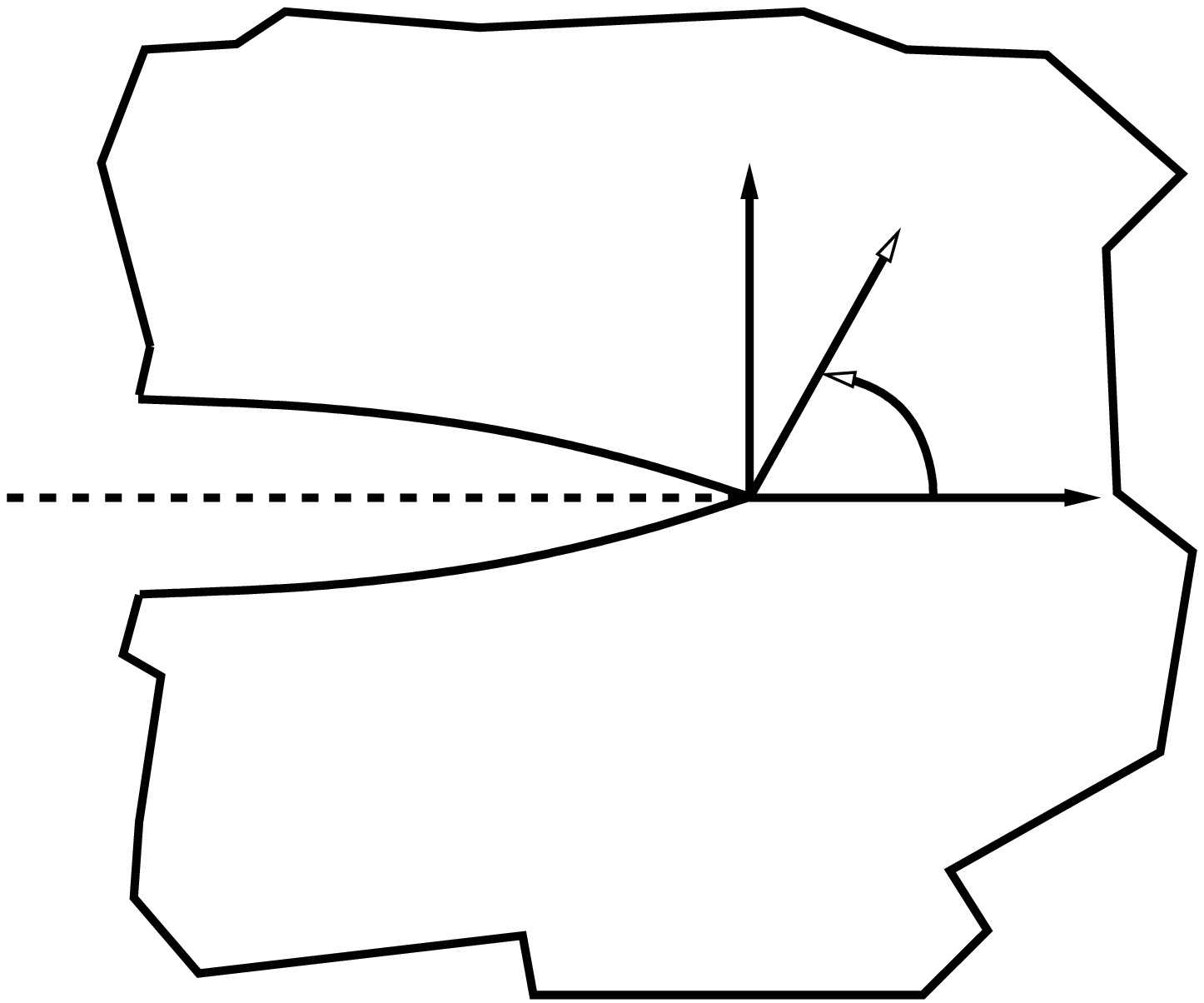}%
\end{picture}%
\setlength{\unitlength}{4144sp}%
\begingroup\makeatletter\ifx\SetFigFont\undefined%
\gdef\SetFigFont#1#2#3#4#5{%
  \reset@font\fontsize{#1}{#2pt}%
  \fontfamily{#3}\fontseries{#4}\fontshape{#5}%
  \selectfont}%
\fi\endgroup%
\begin{picture}(6643,5485)(1022,-6162)
\put(5356,-1726){\makebox(0,0)[lb]{\smash{{\SetFigFont{34}{40.8}{\familydefault}{\mddefault}{\updefault}{\color[rgb]{0,0,0}$x_2$}%
}}}}
\put(5626,-3121){\makebox(0,0)[lb]{\smash{{\SetFigFont{34}{40.8}{\familydefault}{\mddefault}{\updefault}{\color[rgb]{0,0,0}$\varphi$}%
}}}}
\put(6886,-3751){\makebox(0,0)[lb]{\smash{{\SetFigFont{34}{40.8}{\familydefault}{\mddefault}{\updefault}{\color[rgb]{0,0,0}$x_1$}%
}}}}
\put(6121,-2131){\makebox(0,0)[lb]{\smash{{\SetFigFont{34}{40.8}{\familydefault}{\mddefault}{\updefault}{\color[rgb]{0,0,0}$r$}%
}}}}
\end{picture}%